\documentclass[a4paper,leqno,12pt]{article}
\usepackage[english]{babel}
\usepackage[cp1252]{inputenc}
\usepackage[T1]{fontenc}
\usepackage{lmodern}
\usepackage{color}
\usepackage{hyperref}

\usepackage{version}

\usepackage{graphicx}

\excludeversion{pb2}
\excludeversion{pb3}
\excludeversion{notes}
\includeversion{figs} %Figures
\newcommand{\noi}{\noindent }
\newcommand{\noib}{\noindent $\bullet $~}
\newcommand{\pointir}{\discretionary{.}{}{.\kern.35em---\kern.7em}}
\newcommand{\ie}{\emph{i.e. }}

\newcommand{\como}[1]{\bigskip $\blacktriangleright$ \textbf{ \texttt{[#1]}}
\bigskip}

\newcommand{\comf}[1]{\bigskip \textbf{ \texttt{[#1]}} $\blacktriangleleft$ \bigskip}

\usepackage{latexsym}
\usepackage{amsmath}
\usepackage{amssymb}
\newcommand{\Mh}{\widehat{M}}

\newcommand{\gh}{\widehat{g}}
\newcommand{\Dh}{\widehat{D}}

\newcommand{\rich}{\widehat{\mathrm{Ric}}}

\def\buildo#1\over#2{\mathrel{\mathop{\null#2}\limits^{#1}}}
\def\buildu#1\under#2{\mathrel{\mathop{\null#2}\limits_{#1}}}

\newcommand{\R}{\mathbb{R}}
\newcommand{\Rb}{\mathbb{R}^{\bullet}}

\newcommand{\HH}{\mathbb{H}}

\newcommand{\Ss}{\mathbb{S}}

\newcommand{\cC}{{\mathcal C}}
\newcommand{\cD}{{\mathcal D}}

\newcommand{\cK}{{\mathcal K}}

\newcommand{\cS}{{\mathcal S}}

\newcommand{\cU}{{\mathcal U}}

\newcommand{\bra}{\langle}
\newcommand{\ket}{\rangle}

\newcommand{\cty}{C^{\infty}}
\newcommand{\ens}[2]{\{ #1 ~|~ #2 \}}

\newcommand{\pf}{\par{\noindent\textbf{Proof.~}}}

\newbox\qedbox
\setbox\qedbox=\hbox{$\Box$}
\newcommand{\qed}{\hfill\penalty10000\copy\qedbox\par\medskip}
\newtheorem{thm}{Theorem}[section]
\newtheorem{thm-def}[thm]{Theorem and Definition}

\newtheorem{prop}[thm]{Proposition}
\newtheorem{proper}[thm]{Property}

\newtheorem{prop-def}[thm]{Proposition and Definition}
\newtheorem{lem}[thm]{Lemma}

\title{Lindel\"{o}f's theorem for catenoids revisited}
\author{Pierre B\'{e}rard and Ricardo Sa Earp}
\date{July, 2009}

\begin{document}
\maketitle

\pagestyle{myheadings}
%\markboth{B\'{e}rard - Earp}{Work in progress}
%\markright{seulement droite}
\thispagestyle{empty}

\vspace{1.5cm}

\begin{abstract}
\noi In this paper we study the maximal stable domains on minimal
catenoids in Euclidean and hyperbolic spaces and in $\HH^2 \times
\R$. We in particular investigate whether half-vertical catenoids
are maximal stable domains (\emph{Lindel\"{o}f's property}). We also
consider stable domains on catenoid-cousins in hyperbolic space. Our
motivations come from Lindel\"{o}f's 1870 paper on catenoids in Euclidean
space.\\
\end{abstract}

\bigskip

\textbf{Mathematics Subject Classification (2000):~} 53C42, 53C21,
58C40

\medskip

\textbf{Key words:~} Minimal surfaces. Constant mean curvature
surfaces. Stability. Index.

\vspace{2cm}

\label{1-intro}
%%%\section{Introduction}\label{S-intro}
%1-intro

\bigskip
\section{Introduction}\label{S-intro}
\bigskip

\begin{figure}[htbp]
\begin{center}
    \includegraphics[width=6.5cm]{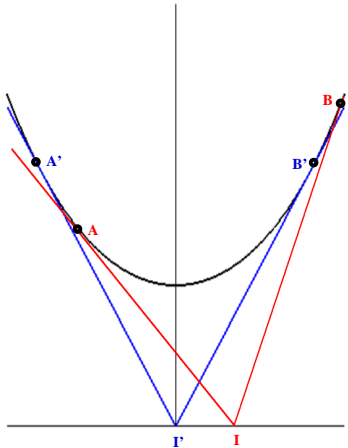}
    \caption[Lindel\"{o}f's construction]{Lindel\"{o}f's construction}
    \label{F-0}
\end{center}
\end{figure}

In his 1870 paper ``Sur les limites entre lesquelles le cat\'{e}no\"{\i}de
est une surface minima'' published in the second volume of the
Mathematische Annalen, see \cite{Lin870}, L. Lindel\"{o}f determines
which domains of revolution on the catenoid $\cC$ in $\R^3$ are
stable. More precisely, he gives the following geometric
construction (see Figure \ref{F-0}).

Take any point $A$ on the generating catenary $C =\ens{(x,z) \in
\R^2}{z= \cosh(x)}$. Draw the tangent to $C$ at the point $A$ and
let $I$ be the intersection point of the tangent with the axis
$\{z=0\}$. From $I$, draw the second tangent to $C$. It touches $C$
at the point $B$. Lindel\"{o}f's result states that the compact
connected arc $AB$ generates a maximal weakly stable domain on the
catenoid $\cC$ in the sense that the second variation of the area
functional for this domain is zero while it is positive for any
smaller domain, and negative for any larger domain. \bigskip

As a consequence, the upper-half of the catenoid, $\cC \cap \{x\ge
0\}$ is a maximal weakly stable among domains invariant under
rotations. We will refer to this latter property as \emph{Lindel\"{o}f's
property}.\bigskip

In this paper, we generalize Lindel\"{o}f's result (not the geometric
construction with tangents), to other catenoid-like surfaces:
catenoids in $\HH^3$ and $\HH^2 \times \R$, embedded
catenoid-cousins (rotation surfaces with constant mean curvature $1$
in $\HH^3(-1)$). \bigskip

The global picture looks as follows. Catenoids in $\R^3$ and
catenoid-cousins in $\HH^3$ satisfy Lindel\"{o}f's property. That
minimal catenoids in $\R^3$ and catenoid-cousins in $\HH^3$ have
similar properties is not surprising from the local correspondence
between minimal surfaces in $\R^3$ and surfaces with constant mean
curvature $1$ in $\HH^3(-1)$. One may observe that the Jacobi
operators look the same, namely $-\Delta - |A_0|^2$, where $A_0$ is
the second fundamental form for catenoids and its traceless analog
for catenoid-cousins. \bigskip

Catenoids in $\HH^2 \times \R$ have index $1$ and do not satisfy
Lindel\"{o}f's property. Catenoids in $\HH^3$ divide into two
families, a family of stable catenoids which foliate the space,
and a family of index $1$ catenoids which intersect each other and
have an envelope. The hyperbolic catenoids do not satisfy
Lindel\"{o}f's property. One may observe that the Jacobi operators
can be written $-\Delta + c - |A|^2$, where $c=1$ for catenoids in
$\HH^2 \times \R$ and $c=2$ for catenoids in $\HH^3$. The presence
of $c$ may explain the extra stability properties of these
catenoids.
\bigskip

To prove his result, Lindel\"{o}f introduced the $1$-parameter family of
Euclidean catenoids passing through a given point and considered the
Jacobi field associated with the variation of this family. In this
paper we work directly with Jacobi fields. More precisely, we
consider the vertical Jacobi field (associated with the translations
along the rotation axis in the ambient space), the variation Jacobi
field (the catenoids come naturally in a $1$-parameter family) and a
linear combination of these two Jacobi fields which is well-suited
to study Lindel\"{o}f's property. \bigskip

In some instances, we could use alternative methods to prove (or
disprove) Lindel\"{o}f's property. For example, the fact that Euclidean
catenoids satisfy Lindel\"{o}f's property follows from the theorem of
Barbosa - do Carmo, see \cite{BC80}, relating the stability of a
domain with the area of its spherical image by the Gauss map. One
could also use the fact that the Jacobi operator on the Euclidean
catenoid is transformed into the operator $-\Delta -2$ on the sphere
minus two points by a conformal map. \bigskip

Such alternative methods are not always available. On the
other-hand, our method applies to catenoids in higher dimensions
as well as to rotation surfaces with constant mean curvature $H$,
$0 \le H \le 1$ in the hyperbolic space $\HH^3(-1)$. These
catenoids or catenoid-like hypersurfaces do not satisfy
Lindel\"{o}f's property. \bigskip

We finally point out that among the examples we have studied, the
hypersurfaces which do not satisfy Lindel\"{o}f's property are
precisely those which are vertically bounded.  \bigskip

Note that the stability of (minimal) catenoids have been studied
in \cite{CD83, Mo81, Seo09} when the ambient space is $\HH^3$, in
\cite{BS08a, ST05} when the ambient space is $\HH^2 \times \R$ and
that the index of catenoid-cousins has been studied in
\cite{LR98}. \bigskip

The paper is organized as follows. In Section \ref{S-pre}, we
recall some basic notations and facts. We review Lindel\"{o}f's
original result in Section \ref{S-r3}. In Section \ref{S-rn}, we
study Lindel\"{o}f's result for $n$-catenoids (minimal rotation
hypersurfaces) in Euclidean space $\R^{n+1}$. In Sections
\ref{S-h2r} and \ref{S-hnr}, we consider catenoids in $\HH^2
\times \R$ and $\HH^n \times \R$. In Section~\ref{S-h3}, we study
minimal catenoids and catenoid-cousins in $\HH^3$.
\bigskip

The authors would like to thank Manfredo do Carmo for pointing out
Lindel\"{o}f's paper to them. \bigskip

The authors would like to thank the Mathematics Department of
PUC-Rio (PB) and the Institut Fourier -- Universit\'{e} Joseph
Fourier (RSA) for their hospitality. They gratefully acknowledge
the financial support of CNPq, FAPERJ, Universit\'{e} Joseph
Fourier and R\'{e}gion Rh\^{o}ne-Alpes.\bigskip

\label{2-pre}
%%%\section{Preliminaries}\label{S-pre}
%2-pre

\bigskip
\section{Preliminaries}\label{S-pre}
\bigskip

Let $(\Mh, \gh)$ be an orientable Riemannian manifold and let $M^n
\looparrowright \Mh^{n+1}$ be a complete orientable minimal
immersion. The second variation of the area functional is given by
the Jacobi operator $J_M$ acting on $\cty_0(M)$,

\begin{equation}\label{E2-1}
J_M = - \Delta - \big( \rich(N,N) + |A|^2\big),
\end{equation}

where $\Delta$ is the (non-positive) Laplacian in the induced metric
on $M$, $N$ a unit normal field along the immersion, $A$ the second
fundamental form of the immersion with respect to $N$ and $\rich$
the Ricci curvature of the ambient space $\Mh$, see \cite{Law80}.
\bigskip

The Jacobi operator appears naturally when one considers families of
constant mean curvature immersions. More precisely, let $X(a,\cdot):
M^n \looparrowright (\widehat{M}^{n+1}, \gh)$ be a $1$-parameter
family of  orientable immersions with unit normal $N(a,x)$ and
constant mean curvature $H(a)$. Let $u(a,x):=\gh(\frac{\partial
X}{\partial a}(a,x), N(a,x))$. Then, see \cite{BGS87},

\begin{equation}\label{E2-2}
H'(a)= \Delta u + \big( \rich(N,N) + |A|^2 \big)u=-J_M(u).
\end{equation} \bigskip

In particular, if $H(a)$ does not depend on $a$, then $u$ satisfies
the equation $J_M(u)=0$. \bigskip

We call \emph{Jacobi field} on $D \subset M$ a $\cty$ function $f$
such that $J_M(f) = 0$ on $D$. The geometry of the ambient space
provides usefull Jacobi fields. More precisely, the following
classical properties follows immediately from Equation (\ref{E2-2}).
\bigskip

\begin{proper}[Killing Jacobi field]\label{P2-1}
Let $M \looparrowright \Mh$ be a minimal or constant mean curvature
immersion and let $\cK$ be a Killing field on $\Mh$. The function
$f_{\cK} = \gh(\cK,N)$, given by the inner product (in $\Mh$) of the
Killing field $\cK$ with the unit normal $N$ to the immersion, is a
Jacobi field on $M$.
\end{proper}\bigskip

\begin{proper}[Variation Jacobi field]\label{P2-2}
Let $X(a,\cdot) : M \looparrowright \Mh$ be a smooth family of
immersions, with the same constant mean curvature $H$, for $a$ in
some interval around $a_0$. Then the function $v =
\gh(\frac{\partial X}{\partial a}(a_0,\cdot),N)$, the scalar product
of the variation vector field of the family with the unit normal $N$
to the immersion $X(a_0,\cdot)$, is a Jacobi field on $X(a_0,M)$.
\end{proper}\bigskip

We say that a domain $D$ on $M$ is \emph{stable} if $\int_M f J_M(f)
\, d\mu_M > 0$ for all $f \in \cty_0 (D)$, where $d\mu_M$ is the
Riemannian measure for the induced metric on $M$. We say that a
domain $D$ on $M$ is \emph{weakly stable} if $\int_M f J_M(f) \,
d\mu_M \ge 0$ for all $f \in \cty_0 (D)$. We say that a relatively
compact open domain $D$ on $M$ has \emph{index} $k$ if the maximal
dimension of subspaces of $\cty_0(D)$ on which $\int_D f J_M(f) \,
d\mu_M$ is negative, is equal to $k$. Finally, we say that an open
domain is maximally weakly stable if it is weakly stable and if any
bigger open domain is not. \bigskip

Let $D \subset M$ be a relatively compact regular open domain, and
let

$$\lambda_1(D) = \inf \ens{\int_{D} f J_M(f) \, d\mu_M}{f
\in \cty_0(D), \int_M f^2 \, d\mu_M =1},$$

be the least eigenvalue of the Jacobi operator $J_M$ with Dirichlet
boundary conditions on $\partial D$. To say that $D$ is weakly
stable but not stable is equivalent to saying that
$\lambda_1(D)=0$.\bigskip

\begin{proper}[Monotonicity]\label{P2-3}
Let $D_1 \subset D_2$ be two relatively compact domains in $M$, such
that $\mathrm{int}(D_2 \setminus D_1) \not = \emptyset$. Then
$\lambda_1(D_1) > \lambda_1(D_2)$. In particular, if $D_2$ is weakly
stable (\ie $\lambda_1(D_2) \ge 0$) and $\mathrm{int}(D_2 \setminus
D_1) \not = \emptyset$, then $D_1$ is stable (\ie $\lambda_1(D_1) >
0$).
\end{proper}\bigskip

\begin{proper}[Stability criterion]\label{P2-4}
A relatively compact domain $D$ is weakly stable if and only if
there exists a positive function $u : D \to \R_+$ such that
$J_M(u)\ge 0$.
\end{proper}\bigskip

Property \ref{P2-3} is the classical monotonicity principle of
Dirichlet eigenvalues. Property \ref{P2-4} follows from the
divergence theorem. \bigskip

\begin{pb2}
\bigskip \como{PB2}

\input{pb2-pre-stability}

\comf{PB2} \bigskip
\end{pb2}

\label{3-r3}
%%%\section{Catenoids in $\R^3$}\label{S-r3}
%3-r3

\bigskip
\section{Catenoids in $\R^3$}\label{S-r3}
\bigskip

We consider the family of catenoids given by the following
parametrization

\begin{equation}\label{E3-1}
X(a,t,\theta) = \big(a \cosh(\frac{t}{a}) \cos \theta , a
\cosh(\frac{t}{a}) \sin \theta , t \big), ~a>0
\end{equation}

and in particular the catenoid $\cC$ given by $X_1$. The unit normal
to $\cC_a$ is

\begin{equation}\label{E3-2}
N(a,t,\theta) = \big( - \dfrac{\cos \theta}{\cosh(\frac{t}{a})} , -
\dfrac{\sin \theta}{\cosh(\frac{t}{a})} , \tanh(\frac{t}{a}) \big).
\end{equation}

\begin{pb2}
\bigskip \como{PB2}

\input{pb2-r3-1}

\comf{PB2}\bigskip
\end{pb2}

The Jacobi operator on $\cC$ is

\begin{equation}\label{E3-3}
J_{\cC} = - \cosh^{-2}(t) \big( \dfrac{\partial^2}{\partial t^2} +
\dfrac{\partial^2}{\partial \theta^2}\big) - 2 \cosh^{-4}(t),
\end{equation}

with \emph{radial part}

\begin{equation}\label{E3-4}
J_{\cC}^0 = - \cosh^{-2}(t) \dfrac{\partial^2}{\partial t^2} - 2
\cosh^{-4}(t) = \cosh^{-2}(t) L_{\cC}.
\end{equation} \bigskip

According to Property \ref{P2-1}, the function

\begin{equation}\label{E3-5}
v(t) = \tanh(t) = \bra \frac{\partial}{\partial z} , N(t,\theta)\ket
\end{equation}

is a Jacobi field on $\cC$. According to Property \ref{P2-2}, the
function

\begin{equation}\label{E3-6}
e(t) = 1 - t \, \tanh (t) = - \bra
\frac{d}{da}X_a|_{a=1}(a,t,\theta) , N(t,\theta)\ket
\end{equation}

is a Jacobi field on $\cC$. \bigskip

\begin{pb2}
\como{PB2}

\input{pb2-r3-2}

\comf{PB2}
\end{pb2}

\begin{thm}\label{T3-1}
Let $\xi_0$ be the positive zero of the function $e(t) = 1 - t \,
\tanh(t)$.
\begin{enumerate}
    \item The domain $\cD_{\xi_0} = X(1,]-\xi_0, \xi_0[, [0,2\pi])$
    is a maximal weakly stable domain on the catenoid $\cC$.
    \item The domain $\cD_+ = X(1,]0, \infty[, [0,2\pi])$ is a
    maximal weakly stable rotation invariant domain on the catenoid $\cC$.
    More precisely, given any $\alpha > 0$, the function
    \begin{equation}\label{E3-7}
    e(\alpha ,t) = v(\alpha) e(t) + e(\alpha) v(t).
    \end{equation}
    has a unique positive zero $\beta(\alpha)$ and the domain
    $\cD_{\alpha,\beta(\alpha)} =
    X(1, ]-\alpha , \beta(\alpha) [, [0,2\pi])$,
    is a maximal weakly stable rotation invariant domain in $\cC$.
    \item The catenoid $\cC$ has index $1$.
\end{enumerate}
\end{thm}\bigskip

\textbf{Proof of Theorem \ref{T3-1}}.\\[4pt]

{1.~} The function $e(t)$ is a Jacobi field on $\cC$ which satisfies
$$\left\{%
\begin{array}{ll}
    J_{\cC}(e) & = 0 \text{~ in~} \cD_{\xi_0},\\
    e|{\partial \cD_{\xi_0}} & = 0. \\
\end{array}%
\right.    $$

It follows that $\lambda_1(\cD_{\xi_0}) = 0$ and hence, any smaller
open domain $\Omega \varsubsetneq \cD_{\xi_0}$ is stable, while any
larger open domain $\cD_{\xi_0} \varsubsetneq \Omega$ is unstable by
Property \ref{P2-3}. \bigskip

{2.~} The function $v(t)$ being positive in the interior of $\cD_+$,
it follows from Property~\ref{P2-4} that $\cD_+$ is weakly stable.
Take any $\alpha > 0$. Because $e$ and $v$ are Jacobi fields, the
function $e(\alpha ,t)$ defined by (\ref{E3-7}) is a Jacobi field
too and satisfies $e(\alpha , -\alpha) = 0$. Because $e(\alpha , \pm
\infty) = - \infty$ and $\frac{\partial e}{\partial t}(\alpha
,-\alpha) \not = 0$, the function $e(\alpha, \cdot)$ must have
another zero $\beta(\alpha) \not = -\alpha$. That this second zero
is unique and positive can be seen directly, or arguing as follows.
First of all, observe that the function $e(\alpha , \cdot)$ cannot
have two negative zeroes or two positive zeroes because $\cD_- =
X(1, ]- \infty , 0[, [0, 2\pi])$ and $\cD_+$ are weakly stable
(equivalently, use the fact that $v \not =0$ for $t\not = 0$ and
Property \ref{P2-4}). The only issue for $e(\alpha, \cdot)$ is to
have one (and only one) negative zero $- \alpha$, and one (and only
one) positive zero $\beta(\alpha)$. It follows that $0$ is the least
eigenvalue of $J_{\cC}$ in $\cD_{\alpha,\beta(\alpha)}$, with
Dirichlet boundary conditions. Hence this domain is maximally stable
(any smaller domain is stable and any larger domain is unstable). It
also follows that $\cD_+$ is a maximal stable domain among rotation
invariant domains.
\bigskip

{3.~} It follows from Assertion 1 (or from Assertion 2) that $\cC$
has index at least one. It also follows from the proof of Assertion
2 that $J_{\cC}^0$, see (\ref{E3-4}), cannot have index bigger than
or equal to $2$. Using the Jacobi field

$$h(t,\theta) = \bra \frac{\partial}{\partial y} , N(t,\theta)\ket =
- \frac{\cos \theta}{\cosh(t)}$$

and Fourier series decomposition in the variable $\theta$, one can
see that the negative eigenvalues of $J_{\cC}$ in rotation invariant
domains can only come from $J_{\cC}^0$, see \cite{BS08a} for a
detailed proof.  \hfill \qed \bigskip

\textbf{Remarks.}
\begin{enumerate}
    \item Using the function $e(\alpha,t)$ defined by (\ref{E3-7}),
    one can recover Lindel\"{o}f's construction, namely that the
    tangents to the catenary $z=\cosh(x)$ at the points $(-\alpha ,
    \cosh(\alpha))$ and $(\beta(\alpha) , \cosh(\beta(\alpha)))$ intersect
    on the axis $\{z=0\}$.
    \item The function $e(\alpha, \cdot)$ can be obtained, up to a
    multiplicative constant, as the Jacobi field arising from the
    variation of the one-parameter family of catenaries passing
    through the given point $(-\alpha , \cosh(\alpha ))$.
    \item A more careful analysis, using for example the fact that
    the catenoid is conformally equivalent to the sphere minus two
    points or the theorem of Barbosa-do Carmo \cite{BC80}, shows that
    the domain $\cD_+$ is a maximal weakly stable
    domain (not only among rotational invariant domains).
\end{enumerate}

\label{3-rn}
%%%\section{Catenoids in $\R^{n+1}$}\label{S-rn}
%3-rn.tex

\bigskip
\section{Catenoids in $\R^{n+1}$}\label{S-rn}
\bigskip

\subsection{The mean curvature equation} \bigskip

We first review the equation of minimal catenoids in $\R^{n+1}$.
\bigskip

Consider the parametrization of a rotation hypersurface about the
axis $\{x_{n+1}\}$ in the Euclidean space $\R^{n+1}$,

\begin{equation}\label{E-rn-2}
\left\{%
\begin{array}{ll}
F &: \R \times S^{n-1} \rightarrow \R^{n+1},\\[6pt]
F &: (t,\omega) \mapsto \big( f(t) \, \omega , t\big), \\
\end{array}%
\right.
\end{equation}

generated by the curve $t \mapsto \big( f(t), t\big)$ in $\R^2_{\{
x_1, x_{n+1}\} }$ (with $f > 0$). In the sequel, we let $f_t$ denote
the derivative of the function $f$ with respect to $t$. \bigskip

\begin{pb2}
\bigskip \como{PB2}

\input{pb2-rn-1}

\comf{PB2} \bigskip
\end{pb2}

The Riemannian metric induced by $F$ is given by

\begin{equation}\label{E-rn-4}
G_F(t,\omega) = \begin{pmatrix}
  1+f_t^2(t) & 0 \\
  0 & f^2(t) \, \mathrm{Id} \\
\end{pmatrix}.
\end{equation} \bigskip

The unit normal to the immersion $F$ is given by

\begin{equation}\label{E-rn-6}
N_F(t,\omega) = (1+f_t^2)^{-1/2} (-\omega , f_t).
\end{equation} \bigskip

\begin{pb2}
\bigskip \como{PB2}

\input{pb2-rn-2}

\comf{PB2} \bigskip
\end{pb2}

We can deduce the equation satisfied by the mean curvature of the
rotation hypersurface parametrized by $F$

\begin{equation}\label{E-rn-10}
n \, H(t) = - f_{tt} (1+f_t^2)^{-3/2} + (n-1)
f^{-1}(1+f_t^2)^{-1/2}.
\end{equation}

In particular, the hypersurface parametrized by $F$ is minimal if
and only if

\begin{equation}\label{E-rn-12}
f \, f_{tt} = (n-1) \, (1+f_t^2)
\end{equation}

(recall that we assume that $f > 0$). A straightforward computation
yields,

\begin{equation}\label{E-rn-14}
\dfrac{d}{dt} \, \big( f^{n-1} \, (1+f_t^2)^{-1/2} \big) = n \, H(t)
\, f^{n-1} \, f_t .
\end{equation} \bigskip

If $F$ is a minimal immersion, \ie if $f$ satisfies (\ref{E-rn-12}),
then $f$ also satisfies

\begin{equation}\label{E-rn-16}
f^{n-1} \, (1+f_t^2)^{-1/2} = C
\end{equation}

for some constant $C$. It follows from (\ref{E-rn-16}) that a
solution $f$ of the differential equation (\ref{E-rn-12}) which does
not vanish at some point never vanishes on its interval of
definition. \bigskip

\begin{lem}\label{L-rn-2}
Let $(I,f)$ be a solution of the differential equation
(\ref{E-rn-12}), where $I$ is some open interval, and $f$ a
function, $f : I \to \R$.
\begin{enumerate}
    \item The pair $(\check{I},\check{f})$, where $\check{I} = \ens{t \in \R}{-t \in
    I}$ and where $\check{f} : \check{I} \to \R$ is defined by
    $\check{f}(t) = f(-t)$, is also a solution of (\ref{E-rn-12}).
    \item The pair $(I_a,f_a)$, where $I_a = \ens{t \in \R}{\frac{t}{a} \in
    I}$ and where $f_a : I_a \to \R$ is defined by $f_a(t) =
    a f(\frac{t}{a})$, is also a solution of (\ref{E-rn-12}).
\end{enumerate}
\end{lem}\bigskip

\pf The proof is straightforward. \qed \bigskip

\textbf{Remark}. The degree-one differential equation can be
obtained directly using the flux formula, see Appendix A. \bigskip

\subsection{Catenoids in $\R^{n+1}$} \bigskip

For $n \ge 2$, let $(I_n,c_n)$ be the maximal solution of the Cauchy
problem

\begin{equation}\label{E-rn-18}
\left\{%
\begin{array}{ll}
f \, f_{tt} =& (n-1) \, (1+f_t^2),\\
f(0) =& 1, \\
f_t(0) =& 0.\\
\end{array}%
\right.
\end{equation}\bigskip

It follows from the first assertion in Lemma \ref{L-rn-2} that the
interval $I_n$ is of the form $I_n = ]-T_n,T_n[$ for some $T_n$ such
that $0 < T_n \le \infty$, and that $t \mapsto c_n(t)$ is an even
smooth function of $t$ which also satisfies

\begin{equation}\label{E-rn-20}
c_n^{n-1}(t) \, \big( 1+c_{n,t}^2(t)\big)^{-1/2} = 1,
\end{equation} \bigskip

where the notation $c_{n,t}$ stands for the derivative of $c_n$ with
respect to $t$. It follows from the above equations that $c_n(t) \ge
1$ on $I_n$, that $c_n$ is strictly increasing on $[0,T_n[$ and that
the limit

$$X_n := \lim_{t\to T_n, \, t < T_n}c_n(t)$$

exists in $\R_+ \cup \{\infty \}$.\bigskip

From (\ref{E-rn-20}) we conclude that

\begin{equation}\label{E-rn-21}
c_{n,t}(t) = \big( c_n^{2n-2}(t) - 1\big)^{1/2}, ~~t \in [0,T_n[.
\end{equation}  \bigskip

Let $d_n(x)$ be the inverse function of the function $t \mapsto
c_n(t)$ from $]0,T_n[$ to $]1,X_n[$, \ie $ d_n \big(c_n(t)\big)
\equiv t$ for $t > 0$. It follows that the derivative $d_{n,x}$
satisfies

$$d_{n,x}(x) = \big( x^{2n-2} - 1\big)^{-1/2}$$

and hence

\begin{equation}\label{E-rn-22}
d_n(x) = \int_1^{x} (u^{2n-2} - 1)^{-1/2} \, du.
\end{equation}

It follows that $X_n = \infty$ and that

\begin{equation}\label{E-rn-26}
T_n = \int_1^{\infty} (u^{2n-2} - 1)^{-1/2} \, du .
\end{equation}\bigskip

Note that $T_2$ is infinite while $T_n$ is finite for $n\ge
3$.\bigskip

By the second assertion of Lemma \ref{L-rn-2}, for $n\ge 2$ and
$a>0$, the maximal solution of the Cauchy problem

\begin{equation}\label{E-rn-18a}
\left\{%
\begin{array}{ll}
f \, f_{tt} =& (n-1) \, (1+f_t^2),\\
f(0) =& a, \\
f_t(0) =& 0,\\
\end{array}%
\right.
\end{equation}

is $\big( ]-aT_n,aT_n[, a c_n(\frac{t}{a})\big)$. \bigskip

We have proved the \bigskip

\begin{prop}\label{P-rn-2}
For $n\ge 2$, the minimal rotation hypersurfaces generated by the
solution curves to Equation (\ref{E-rn-18a}),

\begin{equation}\label{E-rn-30}
F(a,t,\omega) = \big( a c_n(\frac{t}{a})\, \omega  , t \big), ~ a>0,
~ t\in ]-aT_n,aT_n[, ~ \omega \in S^{n-1},
\end{equation}
form a family of minimal catenoids $\cC_a$ in $\R^{n+1}$.
\end{prop}\bigskip

\subsection{Jacobi fields} \bigskip

We consider the minimal immersions (\ref{E-rn-30}). In the next
formulas, we denote the function $c_n$ by $c$ and the value $T_n$ by
$T$ for simplicity. According to (\ref{E-rn-6}), the unit normal to
$\cC_a$ is given by

\begin{equation}\label{E-rn-31}
N(a,t,\omega) = \big( 1 + c_t^2(\frac{t}{a})\big)^{-1/2} \, \big(
-\omega , c_t(\frac{t}{a}) \big).
\end{equation}\bigskip

The \emph{vertical Jacobi field} $v(a,t) := \langle N(a,t,\omega) ,
\frac{\partial}{\partial t}\rangle$ on the catenoid $\cC_a$ satifies

\begin{equation}\label{E-rn-34}
\left\{%
\begin{array}{l}
v(a,t) = v(\frac{t}{a}), ~ \text{where}\\
v \text{ ~ is an odd function},\\
 v(t) = c_t(t) \big( 1 + c_t^2(t)\big)^{-1/2}
= \mathrm{sgn}(t)\big( 1 - c^{2-2n}(t)\big)^{1/2},\\
v(0)=0,\\
\lim_{t \to T_-}v(t) = 1.\\
\end{array}%
\right.
\end{equation} \bigskip

The \emph{variation Jacobi field} $e(a,t) := \langle N(a,t,\omega) ,
\frac{\partial F}{\partial a}\rangle$ on the catenoid $\cC_a$
satisfies

\begin{equation}\label{E-rn-36}
\left\{%
\begin{array}{l}
e(a,t) = e(\frac{t}{a}), ~ \text{where}\\
e \text{ ~ is an even function},\\
e(t) = - c^{2-n}(t) + t v(t),\\
e(0)=-1,\\
\lim_{t \to T_-}e(t) = T.\\
\end{array}%
\right.
\end{equation} \bigskip

Recall that $T = \infty$ when $n=2$ and that $T$ is finite when $n
\ge 3$.

\begin{pb2}
\bigskip \como{PB2}

\input{pb2-rn-3}

\comf{PB2} \bigskip
\end{pb2}

The Jacobi fields $v(t)$ and $e(t)$ satisfy the same Sturm-Liouville
equation on $]0,T[$. Since $v(t)$ does not vanish on $]0,T[$, it
follows from Sturm's intertwining zeroes theorem that $e$ vanishes
once and only once on $]0,T[$.
\bigskip

\begin{prop}\label{P-rn-4}
The family of catenoids in $\R^{n+1}$ admits an envelope which is a
cone whose slope is given by the unique positive zero of the
function $e$.
\end{prop}\bigskip

\begin{figure}[htbp]
%\begin{center}
\begin{minipage}[c]{6.5cm}
    \includegraphics[width=6.5cm]{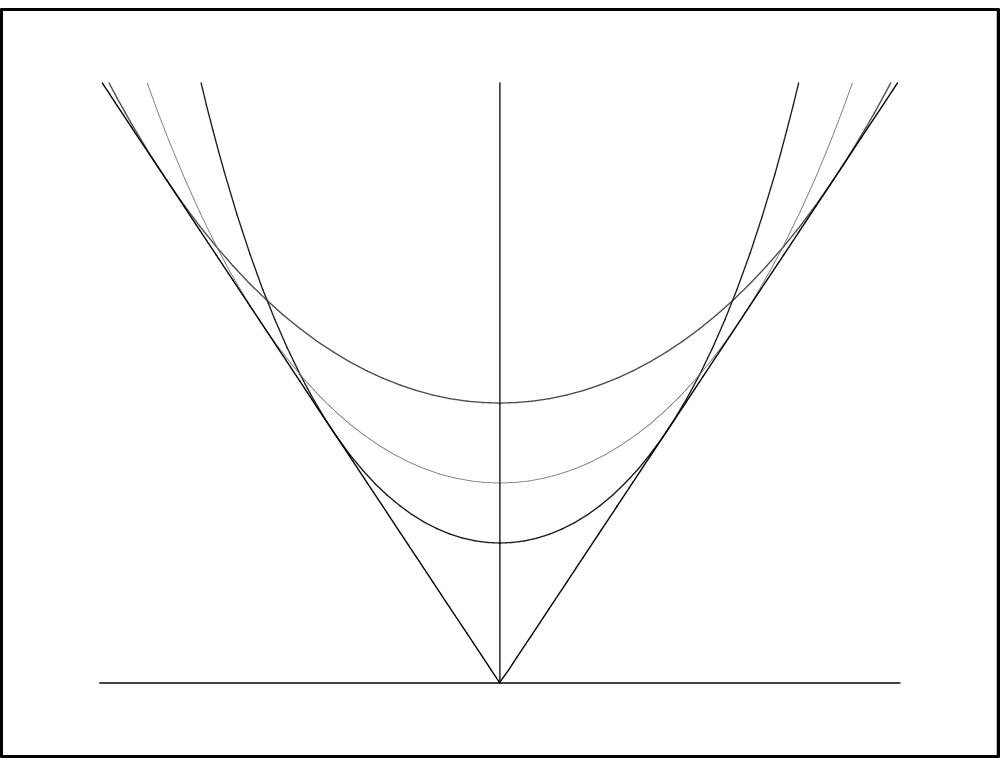}
    \caption[Catenaries and envelope]{Catenaries and envelope}
    \label{F3-rn-env1}
\end{minipage}\hfill
\begin{minipage}[c]{6.5cm}
    \includegraphics[width=6.5cm]{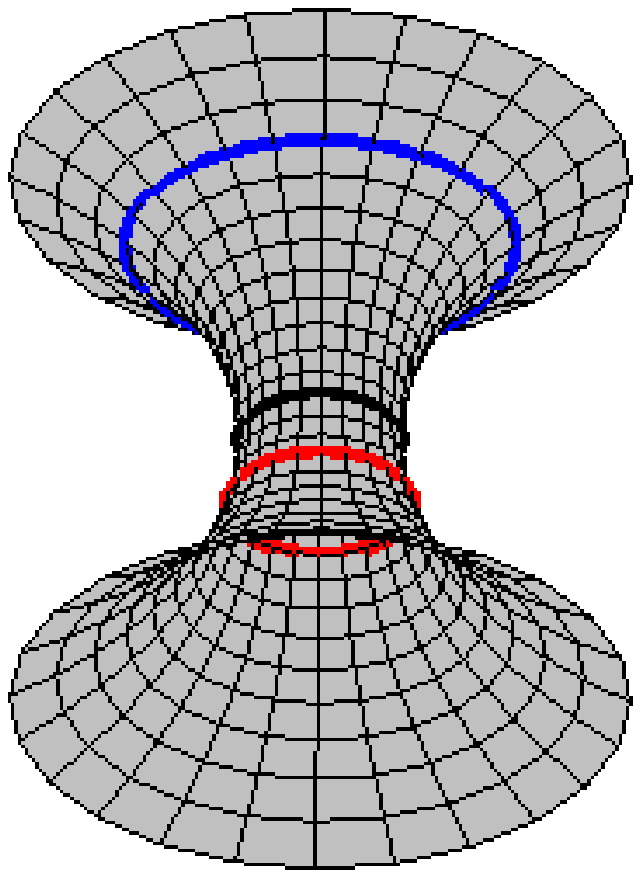}
    \caption[Maximal stability domains]{Maximal stability domains}
    \label{F3-rn-stab1}
\end{minipage}\hfill
%\end{center}
\end{figure}

\pf In order to prove this proposition, it suffices to look at the
family of catenaries which generate the catenoids. The envelope of
this family of curves $\{f(a,t), t \mapsto \big( a
c_n(\frac{t}{a}),t \big)\}_{a>0}$ is given by the equation
$|\frac{\partial f}{\partial a},\frac{\partial f}{\partial t}|=0$,
\ie by the zeroes of the functions $c_n(\frac{t}{a}) - \frac{t}{a}
c_{n,t}(\frac{t}{a})$, \ie by the zeroes $\pm z(a)$ of the functions
$t \mapsto e(a,t)$. Equation (\ref{E-rn-36}) shows that $z(a) = a
z$, where $z$ is the unique positive zero of the Jacobi field $e$.
See Figure \ref{F3-rn-env1}. \qed
\bigskip

 \bigskip

Let $\cC$ denote the catenoid $\cC_1$ and let $F(t,\omega)$ denote
the immersion $F(1,t,\omega)$. \bigskip

\begin{prop}\label{P-rn-6}
The $n$-dimensional catenoid $\cC$ in $\R^{n+1}$ has the following
properties.
\begin{enumerate}
    \item The half catenoid $\cC_+ := F([0,T_n[, S^{n-1})$ is
    weakly stable.
    \item Let $z$ be the unique positive zero of the Jacobi field $e$.
    The domain $\cD_{z} := F(]-z,z[,S^{n-1})$ is a maximal weakly stable
    domain (any larger domain has index at least $1$). It is bounded by the two spheres
    where the catenoid $\cC$ touches the envelope of the family.
    \item The catenoid $\cC$ has index $1$.
\end{enumerate}
\end{prop}\bigskip

\begin{figure}[htbp]
\begin{center}
\begin{minipage}[c]{6.6cm}
    \includegraphics[width=6.6cm]{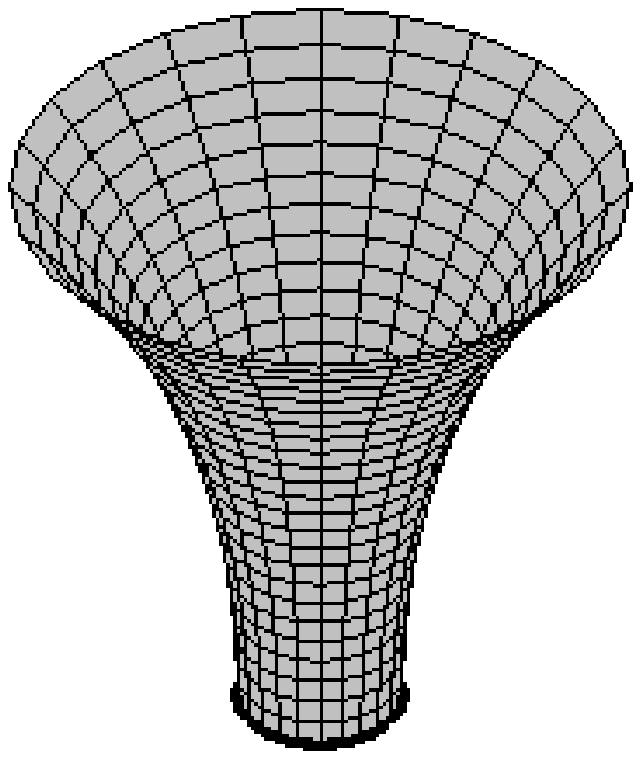}
    \caption[Stable vertical half]{Stable vertical half}
    \label{F3-rn-sv}
\end{minipage}\hfill
\begin{minipage}[c]{6.6cm}
    \includegraphics[width=6.6cm]{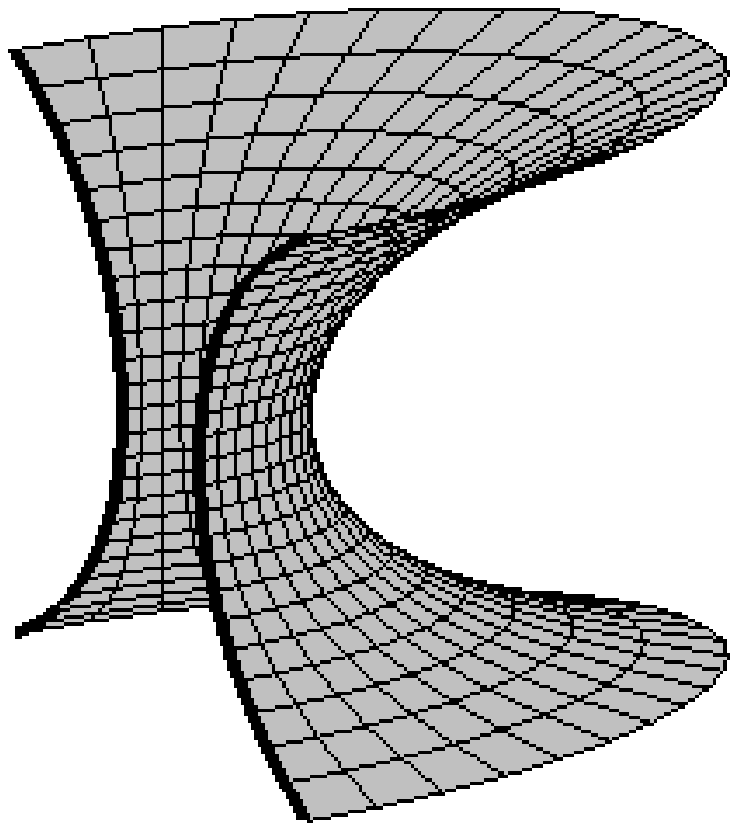}
    \caption[Stable horizontal half]{Stable horizontal half}
    \label{F3-rn-sh}
\end{minipage}\hfill
\end{center}
\end{figure}

\pf Assertions 1 and 2 follow immediately from the properties of the
vertical and variation Jacobi fields and from Properties \ref{P2-3}
and \ref{P2-4}. See Figures \ref{F3-rn-stab1} and \ref{F3-rn-sv}.
\bigskip

Assertion 3 has been proved in \cite{BS08a} using the fact that the
horizontal half-catenoids are stable, see Figure \ref{F3-rn-sh}, and
in \cite{TZ07} by another method. \qed
\bigskip

\begin{thm}\label{T-rn-8}
For $n \ge 3$, the $n$-dimensional catenoid $\cC$ in $\R^{n+1}$ does
not satisfy Lindel\"{o}f's property. More precisely, letting $z$ be the
positive zero of the Jacobi field $e$, there exists an $\ell \in
]0,z[$ such that the following properties hold.
    \begin{enumerate}
    \item The domain $\cD'_{\ell}:= F(]- \ell, \infty[,S^{n-1})$ is weakly
    stable.
    \item For any $\alpha > \ell$, there exists $\beta(\alpha)
    \in ]0,\infty[$ such that the domain $\cD_{\alpha, \beta(\alpha)} :=
    F(]-\alpha,\beta(\alpha)[,S^{n-1})$ is a maximal weakly
    stable domain. In particular, for $\alpha > \ell$, the domain
    $\cD'_{\alpha} := F(]- \alpha, \infty[,S^{n-1})$ has index $1$.
    \item When it exists, the maximal weakly stable domain $\cD_{\alpha,
    \beta(\alpha)}$ is given by Lindel\"{o}f's construction. More
    precisely, the tangents to the catenary $t \mapsto \big( c(t),t \big)$
    at the two points $\big( c(\alpha), -\alpha \big)$ and
    $\big( c(\beta(\alpha)), \beta(\alpha) \big)$ meet on the axis of the catenary.
    \end{enumerate}
\end{thm}\bigskip

\pf

We introduce the Jacobi field

$$w(\alpha,t) := v(\alpha) e(t) + e(\alpha) v(t).$$

Because $e$ is even and $v$ odd, it follows that $w(\alpha,
-\alpha)=0$. Since $w(\alpha,0) = -v(\alpha) < 0$ and $\lim_{t \to
T_n, \, t < T_n}w(\alpha, t) = e(\alpha) + T_n v(\alpha)$ (use
(\ref{E-rn-34}) and (\ref{E-rn-36})). Since $n\ge 3$, the value
$T_n$ is finite and we introduce the Jacobi field $y(t) = e(t) + T_n
v(t)$. We have $y(0)=-1$ and $y(z) = T_n v(z)$. It follows that $y$
has one (and only one) zero $\ell \in ]0,z[$. For $\alpha \le \ell$,
we have that $y(\alpha) \le 0$ and we may conclude that $w(\alpha,
\cdot)$ does not vanish on $]-\alpha , \infty[$. On the other-hand,
when $\alpha > \ell$, $w(\alpha, \cdot)$ has a unique positive zero
$\beta(\alpha)$. \bigskip

3) Writing the equations for the tangents to the catenary at the
points $\big( c(\alpha), -\alpha \big)$ and $\big( c(\beta), \beta
\big)$, we see that a necessary and sufficient condition for the
tangents to intersect on the axis of the catenary is

$$\alpha + \beta = \frac{c(\alpha)}{c_t(\alpha)} +
\frac{c(\beta)}{c_t(\beta)}.$$

Writing that $w(\alpha, \beta)=0$ we find the same necessary and
sufficient condition. This proves the last assertion. \qed
\bigskip

\begin{figure}[htbp]
\begin{center}
\begin{minipage}[c]{6.5cm}
    \includegraphics[width=6cm]{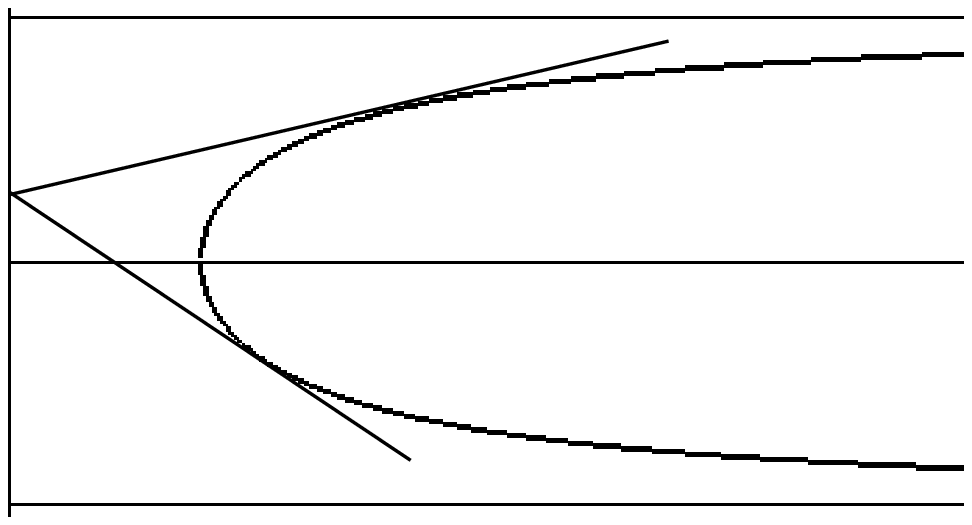}
    \caption[$\cD_{\alpha,\beta(\alpha)}, ~n \ge 3$]{$\cD_{\alpha,
    \beta(\alpha)}, ~n \ge 3$}
    \label{F-cat-Rn-1}
\end{minipage}\hspace{1cm}
\begin{minipage}[c]{6.5cm}
    \includegraphics[width=6cm]{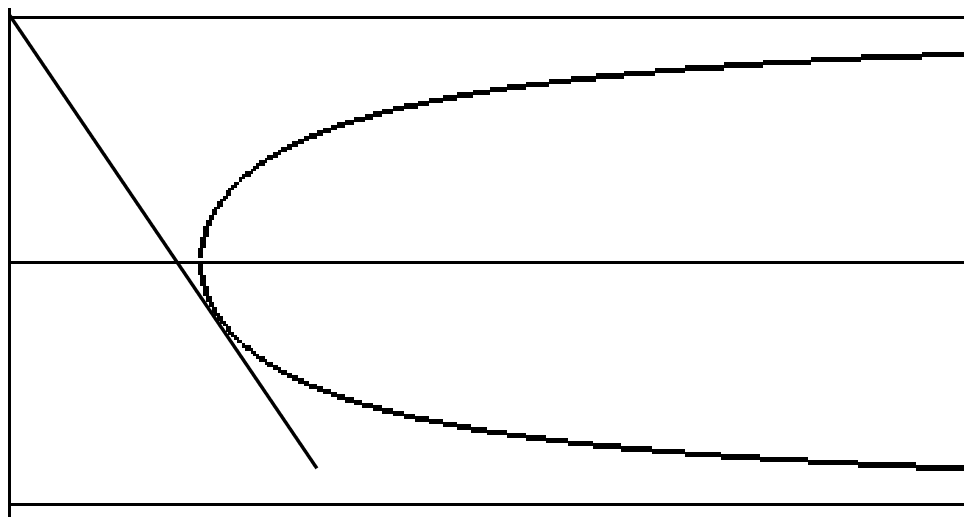}
    \caption[$\cD'_{\ell}, ~n \ge 3$]{$\cD'_{\ell}, ~n \ge 3$}
    \label{F-cat-Rn-2}
\end{minipage}\hfill
\end{center}
\end{figure}

The following figures illustrate the difference between the case
$n \ge 3$ (Figures \ref{F-cat-Rn-1} and \ref{F-cat-Rn-2}) and the
case $n=2$. \bigskip

When $n\ge 3$, the construction of the maximal weakly stable domains
with tangents and the fact that the height of the catenoids is
bounded shows that Lindel\"{o}f's property does not hold.\bigskip

%\newpage

\label{4-h2r}
%%%\section{Catenoids in $\HH^2 \times \R$}\label{S-h2r}
%4-h2r

\bigskip
\section{Catenoids in $\HH^2 \times \R$}\label{S-h2r}
\bigskip

Catenoids in $\HH^2 \times \R$ have been studied in \cite{ST05,
BS08a}. We take the ball model for $\HH^2$ and we let $\rho$ denote
the hyperbolic distance to $0$. We equip $\Mh = \HH^2 \times \R$
with the product metric $\gh$. \bigskip

\subsection{Preliminaries}

In this Section, we review the computations of \cite{BS08a}.
\bigskip

The catenoids in $\HH^2 \times \R$ are generated by catenaries
$C_1(a,\rho) =\big( \tanh(\rho) , \lambda (a,\rho)\big)$ in a
vertical plane $\gamma \times \R$, where $\gamma$ is a complete
geodesic in $\HH^2$, and where

\begin{equation}\label{E4-1}
\lambda(a, \rho) = \sinh(a) \, \int_a^{\rho} \big( \sinh^2(t) -
\sinh^2(a) \big)^{-1/2} \, dt, ~ a>0.
\end{equation}\bigskip

As a matter of fact, $C_1(a,\rho)$ describes a half-catenary and the
whole catenary can be parametrized in the arc-length parameter $s$,
by

\begin{equation}\label{E4-2}
C_2(a,s) = \big( \tanh(R(a,s)/2) , \Lambda(a,s) \big)
\end{equation}

where the function $R(a,s)$ and $\Lambda(a,s)$ are smooth and,
respectively even and odd. They satisfy the relations

\begin{equation}\label{E4-3}
\left\{%
\begin{array}{ll}
R(a,s)& = a + \cosh(a) \int_0^s \sinh(t) \big( \cosh^2(a) \cosh^2(t)
- 1 \big)^{-1/2} \, dt, \\[7pt]
\cosh(R(a,s)) & = \cosh(a) \cosh(s) \text{ ~and~ } R(a,s) \ge a,\\[7pt]
\Lambda(a,s) & = \sinh(a) \displaystyle \int_0^s
(\cosh^2(a) \cosh^2(t) -1)^{-1/2} \, dt,\\[7pt]
\Lambda_s^2 + R_s^2 & \equiv 1.\\
\end{array}%
\right.
\end{equation}\bigskip

The family $\{\cC_a\}_{a>0}$ of catenoids in $\HH^2 \times \R$ is
given (in the ball model) by

\begin{equation}\label{E4-4}
X(a,s,\theta) =
\begin{pmatrix}
\tanh(R(a,s)/2) \omega_{\theta}\\[6pt]
\Lambda(a,s)\\
\end{pmatrix}.
\end{equation}

The metric $X^{*}\gh$ induced by $X$ on $\cC_a$ is $ds^2 +
\sinh^2(R(a,s)) d\theta^2$ and the unit normal is given by

\begin{equation}\label{E4-5}
N(a,s,\theta) =
\begin{pmatrix}
\dfrac{\Lambda_s}{2 \cosh^2(R(a,s)/2)} \omega_{\theta}\\
\\
R_s\\
\end{pmatrix},
\end{equation}

where $\Lambda_s = \frac{\partial \Lambda}{\partial s}$ and $R_s =
\frac{\partial R}{\partial s}$. \bigskip

\subsection{Jacobi fields}\bigskip

\noib The \emph{vertical Jacobi field} is the function

\begin{equation}\label{E4-6}
v(a,s) = \gh(\frac{\partial}{\partial t},N) = R_s.
\end{equation}

Taking (\ref{E4-3}) into account, we find

\begin{equation}\label{E4-7}
v(a,s) = \cosh(a) \sinh(s) \big( \cosh^2(a) \cosh^2(s) - 1
\big)^{-1/2}.
\end{equation}

Note that $v(a,0) = 0$ and $v(a, \infty) = 1$. \bigskip

\noib We take the \emph{variation Jacobi field} to be

\begin{equation}\label{E4-8}
e(a,s) = - \gh(\frac{\partial X}{\partial a}(a,s,\theta),
N(a,s,\theta)).
\end{equation}

This Jacobi field has been computed in \cite{BS08a}. We have

$$e(a,s) = \Lambda_s R_a - \Lambda_a R_s, \text{~~ and}$$

\begin{equation}\label{E4-9}
\begin{array}{ll}
e(a,s) = & \sinh^2(a) \cosh(s) \big( \cosh^2(a) \cosh^2(s) - 1
\big)^{-1} \cdots \\
    & \hphantom{xxxxx} - v(a,s) \displaystyle \int_0^s
    B(a,t) \, dt,
\end{array}
\end{equation}

where $$B(a,t) := \cosh(a)  \sinh^2(t) \big( \cosh^2(a) \cosh^2(t) -
1 \big)^{-3/2}.$$ \bigskip

\subsection{Stable domains on $\cC_a$}\bigskip

Define the rotation invariant domains

\begin{equation}\label{E4-10}
\cD_{\pm} = X(a, \R_{\pm},[0, 2\pi]), ~\text{and}
\end{equation}

\begin{equation}\label{E4-11}
\cD_{\alpha} = X(a, ]-\alpha ,\alpha [,[0, 2\pi]).
\end{equation}

In \cite{BS08a}, we proved the following result.

\begin{thm}\label{T4-1}
The catenaries $\cC_a \subset \HH^2 \times \R$ have the following
properties.
\begin{enumerate}
    \item The domains $\cD_{\pm}$ are weakly stable.
    \item The function $e(a,s)$ has a unique positive zero $z(a)$,
    and
    \begin{itemize}
        \item $\cD_{\alpha}$ is stable for $0 < \alpha  < z(a)$.
        \item $J_{\cC_a}$ has eigenvalue $0$ in $\cD_{z(a)}$ with
        Dirichlet boundary conditions.
        \item $\cD_{\alpha}$ is unstable for $\alpha > z(a)$.
    \end{itemize}
    \item For all $a > 0$, the catenoid $\cC_a$ has index $1$.
\end{enumerate}
\end{thm}\bigskip

\textbf{Sketch of the proof of Theorem \ref{T4-1}}.

\emph{Assertion 1} follows from Property \ref{P2-4}, using the
Jacobi field $v(a,s)$ which does not vanish in the interior of
$\cD_{\pm}$.

\emph{Assertion 2} is a consequence of Property \ref{P2-1} and the
fact that the function $e(a,s)$ has a (unique) zero on $]0, +
\infty[$. Note that the uniqueness of the positive zero of $e(a,s)$
is a consequence of Assertion 1.

\emph{Assertion 3}. We refer to \cite{BS08a}. \qed \bigskip

\begin{thm}\label{T4-2}
The catenoids $\cC_a$ in $\HH^2 \times \R$ do not satisfy Lindel\"{o}f's
property: the domains $\cD_{\pm}$ are not maximally weakly stable. More
precisely, there exists a unique $\ell(a) \in ]0,z(a)[$ such that
$\cD'_{\ell(a)} := X(a,]-\ell(a), \infty[,[0,2\pi])$ is maximally
weakly stable among rotationally invariant domains.
\end{thm} \bigskip

\pf For $\alpha > 0$, introduce the Jacobi field

\begin{equation}\label{E4-12}
e(a,\alpha,s) = v(a,\alpha) e(a,s) + e(a,\alpha) v(a,s)
\end{equation}

where $e$ and $v$ are given by (\ref{E4-6}) and (\ref{E4-9}).
\bigskip

Because the vertical Jacobi field $v$ does not vanish on
$]0,\infty[$ and on $]-\infty,0[$, we have that $e(a,\alpha, \cdot)$
has at most one zero on these intervals (see also Theorem
\ref{T4-1}, Assertion 1). Oberve that $e(a,\alpha,-\alpha)=0$ and
that $e(a,\alpha,0)=v(a,\alpha)>0$. It follows that $e(a,\alpha,
\cdot)$ has a zero in $]0,\infty[$ if and only if $e(a,\alpha,
\cdot)$ is negative near infinity. Using (\ref{E4-12}), we can write
$e(a,s) = f(a,s) - v(a,s) \int_0^s B(a,t) \, dt$ where $f(a,s) =
\sinh^2(a) \cosh(s) \big( \cosh^2(a) \cosh^2(s) - 1 \big)^{-1}$,
$f(a,0)=1$, $f(a,\infty)=0$ and $B(a,t) = \cosh(a) \sinh^2(t) \big(
\cosh^2(a) \cosh^2(t) - 1 \big)^{-3/2}$. Let $E(a) = \int_0^{\infty}
B(a,t) \, dt$, a positive finite value. Using these notations, we
have

$$e(a,\alpha,s) = v(a,\alpha) f(a,s) + v(a,s) [e(a,\alpha) -
v(a,\alpha) \int_0^s B(a,t) \, dt].$$

The sign of $e(a,\alpha,s)$ near $+\infty$ is given by the sign of
$e(a,\alpha) - E(a) v(a,\alpha)$. \bigskip

~~$\rhd$~ If $\alpha > z(a)$ (the unique positive zero of the
variation Jacobi field $e(a,\cdot)$, then $e(a,\alpha) < 0$ so that
$e(a,\alpha) - E(a) v(a,\alpha) < 0$ and hence $e(a,\alpha,s)$ must
have a zero $\beta(\alpha) \in ]0,\infty[$. Clearly, we must have $0
< \beta(\alpha) < z(a)$. This is not surprising in view of Theorem
\ref{T4-1}, Assertion 2. \bigskip

~~$\rhd$~ If $\alpha = z(a)$, then $e(a,z(a),s) = v(a,z(a))e(a,s)$
has two zeroes $\pm z(a)$. \bigskip

~~$\rhd$~ If $0 < \alpha < z(a)$, consider the Jacobi field $w(t):=
e(a,t) - E(a) v(a,t)$. We have $w(0)=1$ and $w(z(a))=-E(a) v(a,
z(a)) < 0$ so that $w$ has a unique positive zero $\ell(a) \in
]0,z(a)[$. When $0 < \alpha < \ell(a), w(\alpha ) > 0$ and hence
$X(a, ]-\alpha ,\infty[,[0,2\pi]$ is weakly stable. For $\alpha >
\ell(a), w(a)<0$, $e(a,\alpha,\cdot)$ has a positive zero $\beta
(\alpha)$ and $X(a, ]-\alpha ,\beta (\alpha)[,[0,2\pi]$ is a maximal
weakly stable domain. \qed \bigskip

\textbf{Remark.} One could show that $\cD_{\ell(a)}$ is maximally
weakly stable by using a conformal transformation. This method does
not apply in higher dimension whereas the above one does. \bigskip

\label{4-hnr}
%%%\section{Catenoids in $\HH^n \times \R$}\label{S-hnr}
%4-hnr

\bigskip
\section{Catenoids in $\HH^n \times \R$}\label{S-hnr}
\bigskip

We consider the space $\HH^n \times \R$ with the product metric $\gh
= g_h + dt^2$ and we work with the ball model for $(\HH^n, g_h)$.
\bigskip

We consider a rotation hypersurface about the axis $\R$, with
parametrization

\begin{equation}\label{E-hnr-2}
F(t,\omega) = \big( \tanh(f(t)/2)\omega , t \big)
\end{equation}

where $f(t) > 0$ is the hyperbolic distance to the axis. Using the
flux formula (see Appendix A), we obtain easily the following
differential equation for minimal rotation hypersurfaces in $\HH^n
\times \R$,

\begin{equation}\label{E-hrn-4}
\sinh^{n-1}\big( f(t) \big) \big( 1 + f_t^2(t) \big)^{-1/2} = C
\end{equation}

for some constant $C$, where $f_t$ denotes the derivative of $f$
with respect to $t$. \bigskip

Differentiating this equation, we have that $f$ also satisfies the
equation

\begin{equation}\label{E-hrn-6}
\sinh\big( f(t) \big) f_{tt}(t) - (n-1) \cosh\big( f(t) \big) \big(
1+f_t^2\big) = 0.
\end{equation}

\begin{lem}\label{L-hnr-2}
The Cauchy problem
\begin{equation}\label{E-hrn-8}
\left\{%
\begin{array}{lll}
\sinh\big( f(t) \big) f_{tt}(t)& = & (n-1) \cosh\big( f(t) \big)
\big( 1+f_t^2\big),\\
    f(0) & = & a,\\
    f_t(0) & = & 0,\\
\end{array}%
\right.
\end{equation}
has a maximal solution of the form $\big( ]-T(a),T(a)[, f(a,t)
\big)$ where $t \mapsto f(a,t)$ is a smooth, even function of $t$.
Furthermore, the function $f$ satisfies

\begin{equation}\label{E-hrn-10}
\left\{%
\begin{array}{lll}
\sinh^{n-1}\big( f(a,t)\big) \big( 1+f_t^2\big)^{-1/2} & = & \sinh^{n-1}(a),\\
    f(a,t) & \ge & a, \text{ ~for all~ } t,\\
    f_t  \text{ ~~has the sign of~ } t.\\
\end{array}%
\right.
\end{equation}
\end{lem}\bigskip

For $t \ge 0$, we have

$$f_t(a,t) = \big( \dfrac{\sinh^{2n-2}(f(a,t)) - \sinh^{2n-2}(a)}{\sinh^{2n-2}(a)}
\big)^{1/2}$$

and the function $f(a,\cdot)$ is a bijection from $[0,T(a)[$ to $[a,
\infty[$. Let $\lambda(a, \cdot) : [a, \infty[ \to [0,T(a)[$ be the
inverse function to $f$. Then

\begin{equation}\label{E-hrn-12}
\lambda(a,\rho) = \sinh^{n-1}(a) \int_a^{\rho} \big( \sinh^{2n-2}(u)
- \sinh^{2n-2}(a)\big)^{-1/2} \, du
\end{equation}

which shows that the value $T(a)$ is finite,

\begin{equation}\label{E-hrn-14}
T(a) = \sinh^{n-1}(a) \int_a^{\infty} \big( \sinh^{2n-2}(u) -
\sinh^{2n-2}(a)\big)^{-1/2} \, du
\end{equation}\bigskip

It follows from the preceding formulas that

\begin{equation}\label{E-hrn-16}
f_t \big( 1 + f_t^2 \big)^{-1/2} = \mathrm{sgn}(t) \big( 1 -
(\frac{\sinh(a)}{\sinh(f)})^{2n-2} \big)^{1/2}.
\end{equation} \bigskip

\textbf{Jacobi fields} \bigskip

The normal to the catenoid is given by

\begin{equation}\label{E-hrn-18}
N(a,t,\omega) = (1+f_t^2)^{-1/2} \big( - \frac{\omega}{2
\cosh^2(f/2)},f_t \big).
\end{equation}

It follows that the \emph{vertical Jacobi field} $v(a,t)$ is an odd
function of $t$ which satisfies

\begin{equation}\label{E-hrn-20}
\left\{%
\begin{array}{ll}
v(a,t) & := \gh(N,\frac{\partial}{\partial t}) = f_t \big( 1 + f_t^2
\big)^{-1/2}\\[6pt]
& = \mathrm{sgn}(t) \big( 1 - (\frac{\sinh(a)}{\sinh(f)})^{2n-2}
\big)^{1/2}, \\[6pt]
v(a,0) &= 0, ~~\lim_{t \to T(a)} v(a,t)=1.\\
\end{array}%
\right.
\end{equation}\bigskip

The \emph{variation Jacobi field} $e(a,t)$ is an even function of
$t$ which satisfies

\begin{equation}\label{E-hrn-22}
\left\{%
\begin{array}{ll}
e(a,t) & := \gh(N,\frac{\partial F}{\partial a}) = - f_a \big( 1 +
f_t^2 \big)^{-1/2}\\[6pt]
& = f_t (1+f_t^2)^{-1/2}\lambda_a(a,f), \text{ ~for~ } t\ge 0,\\[6pt]
v(a,0) &= -1.\\
\end{array}%
\right.
\end{equation}\bigskip

Note that the second equality follows from the fact that
$f(a,\lambda(a,\rho)) \equiv \rho$ for $\rho > a$. \bigskip

It follows from Equation (\ref{E-hrn-12}) that

$$
\begin{array}{ll}
\lambda_a(a,\rho) &= - \; \frac{\sinh(\rho)}{\sinh(a)}
\frac{\cosh(a)}{\cosh(\rho)} \big(
(\frac{\sinh(\rho)}{\sinh(a)})^{2n-2} - 1 \big)^{-1/2} + ~~\cdots\\[8pt]
& \hphantom{xxxx} + \cosh(a) \int_1^{\sinh(\rho)/\sinh(a)}
(v^{2n-2}-1)^{-1/2} (\sinh^2(a)v^2+1)^{-3/2} \, dv.
\end{array}
$$

Using the above expressions, we find that for $t \ge 0$,

\begin{equation*}\label{E-hrn-24}
\begin{array}{ll}
e(a,t) & = - \; \frac{\cosh(a)}{\cosh(f)} \big(
\frac{\sinh(a)}{\sinh(f)} \big)^{n-2} + ~~ \cdots\\[6pt]
& \hphantom{xxx} v(a,t) \cosh(a) \int_1^{\sinh(f)/\sinh(a)}
(v^{2n-2}-1)^{-1/2} (\sinh^2(a)v^2+1)^{-3/2} \, dv.
\end{array}
\end{equation*}\bigskip

We write the preceding equality as

\begin{equation}\label{E-hrn-26}
\left\{%
\begin{array}{ll}
e(a,t) & =: - e_0(a,t) + v(a,t) e_1(a,t), \text{ ~where,}\\[6pt]
e_0(a,t) &:= \frac{\cosh(a)}{\cosh(f)} \big(
\frac{\sinh(a)}{\sinh(f)} \big)^{n-2}, \text{ ~positive and even,}\\[6pt]
e_0(a,0) & = 1, ~~e_0(a,T(a)-)=0,\\[6pt]
e_1(a,0) & = \cosh(a) \int_1^{\sinh(f)/\sinh(a)}
(v^{2n-2}-1)^{-1/2} (\sinh^2(a)v^2+1)^{-3/2} \, dv\\[6pt]
e_1(a,0) &= 0, ~~e_1(a,T(a)-) = E(a),\\[6pt]
\end{array}%
\right.
\end{equation}

where $E(a) := \cosh(a) \int_1^{\infty} (v^{2n-2}-1)^{-1/2}
(\sinh^2(a)v^2+1)^{-3/2} \, dv$ is a finite, positive value.
\bigskip

\begin{prop}\label{P-hnr-1}
The vertical Jacobi field $v(a,t)$ only vanishes at $t=0$. As a
consequence, the half-vertical catenoids $\cC_{a,\pm} :=
F(a,\R_{\pm},S^{n-1})$ are weakly stable.\\[6pt]
The variation Jacobi field $e(a,t)$ has exactly one positive zero
$z(a) \in ]0,T(a)[$. As a consequence the domain $\cD_{z(a)} :=
F(a,]-z(a),z(a)[,S^{n-1})$ is a maximal weakly stable domain.
\end{prop}\bigskip

\pf The proof is clear in view of Properties \ref{P2-3} and
\ref{P2-4}. Note that the fact that $e(a,\cdot)$ has a unique
positive zero follows from the positivity of $v(a,\cdot)$ in
$]0,\infty[$ and Sturm intertwining zeroes theorem. \qed \bigskip

We now introduce the Jacobi field

$$e(a,\alpha,t) := v(a,\alpha) e(a,t) + e(a,\alpha) v(a,t).$$

Notice that $e(a,\alpha, -\alpha)=0$ and that $e(a,\alpha,0) =
-v(a,\alpha) < 0$, so that $e(a,\alpha,\cdot)$ cannot have another
zero on $]-\infty,0[$. For $t \ge 0$, consider the Jacobi field

$$y(a,t):= e(a,t) + E(a) v(a,t).$$

We have that $$y(a,\alpha) = \lim_{t \to T(a)-} e(a,\alpha,t).$$

It is clear that $y(a,\cdot)$ has a unique zero on $]0, \infty[$,
namely some $\ell(a) \in ]0,z(a)[$ (where $z(a)$ is the positive
zero of $e(a,\cdot)$). \bigskip

For $0 < \alpha  \le \ell(a)$, we have that $y(a,\alpha) \le 0$ and
hence that $e(a,\alpha,t) < 0$ for $t$ close enough to $T(a)$. This
implies that for such values of $\alpha$, the function
$v(a,\alpha,\cdot)$ cannot vanish on $]0,T(a)[$. For $\alpha >
\ell(a)$, we have that $y(a,\alpha) > 0$ so that $e(a,\alpha,
\cdot)$ has a (unique) zero $\beta(\alpha) \in ]0,T(a)[$. \bigskip

We have proved, \bigskip

\begin{prop}\label{P-hnr-2}
With the above notations,
\begin{enumerate}
    \item the domain $\cD'_{\ell(a)} := F(a, ]-l(a),T(a)[,S^{n-1})$
    is a maximal rotationally symmetric weakly stable domain,
    \item the domain $\cD_{\alpha,\beta(\alpha)} := F(a,]-\alpha,
    \beta(\alpha)[,S^{n-1})$ is a maximal weakly stable domain.
\end{enumerate}
\end{prop}\bigskip

\label{5-h3}
%%%\section{Catenoids and catenoid cousins in $\HH^3$}\label{S-h3}
%5-h3.tex

\bigskip
\section{Catenoids and catenoid cousins in $\HH^3$}\label{S-h3}
\bigskip

\subsection{Hyperbolic computations}
\bigskip

We work in the half-space model for the hyperbolic space,
\begin{equation}\label{E5-1-1}
\HH^3_{\{x_1,x_2,x_3\}} = \ens{(x_1,x_2,x_3) \in \R^3}{x_3 > 0},
\hspace{0.5cm} g_h = x_3^{-2}\big( dx_1^2 + dx_2^2 + dx_3^2 \big).
\end{equation}\bigskip

\begin{figure}[h!tb]
\begin{center}
    \includegraphics[scale=0.5, angle=0]{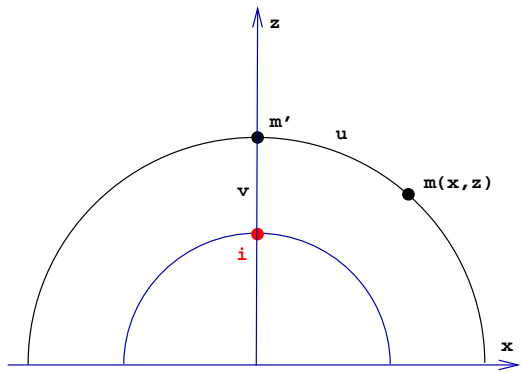}
    \caption[Hyperbolic $2$-plane]{Hyperbolic $2$-plane}
    \label{F5-1-1}
\end{center}
\end{figure}

In the hyperbolic plane

\begin{equation}\label{E5-1-2}
\HH^2_{\{x,z\}} = \ens{(x,z) \in \R^2}{z>0}, \hspace{0.5cm} g_h =
z^{-2}\big( dx^2 + dz^2\big),
\end{equation}

we consider the Fermi coordinates $(u,v)$ defined as follows (see
Figure \ref{F5-1-1}). Given a point $m=m(x,z)$, let $m'$ be its
orthogonal projection on the vertical geodesic $\gamma =
\ens{(0,e^t)}{t \in \R} \subset \HH^2_{\{x,z\}}$. Let $u$ denote the
signed hyperbolic distance $d_h(m,m')$ and $v$ the signed hyperbolic
distance $d_h(m',i)$ (where $i=(0,1)$).\bigskip

\bigskip

The following formulas relate the coordinates $(x,z)$ to the
coordinates $(u,v)$.

\begin{equation}\label{E5-1-3}
\begin{array}{lll}
\left\{%
\begin{array}{lll}
x &=& e^v \tanh(u),\\
z &=& \dfrac{e^v}{\cosh(u)},\\
\end{array}%
\right. &\text{ ~and~ }&
\left\{%
\begin{array}{lll}
u &=& \arg \sinh (\dfrac{x}{z}),\\
v &=& \dfrac{1}{2}\ln (x^2+z^2).\\
\end{array}%
\right. \\
\end{array}\bigskip
\end{equation}\bigskip

In the coordinates $\{u,v\}$, the hyperbolic metric is given by

\begin{equation}\label{E5-1-4}
g_h = du^2 + \cosh^2(u) dv^2.
\end{equation}\bigskip

\bigskip
\subsection{Rotation surfaces in $\HH^3$}
\bigskip

We consider a curve $f(t)=(t,f(t))$ in the plane $\HH^2_{\{u,v\}}$
and the corresponding rotation surface $F : M
\looparrowright\HH^3_{\{x_1,x_2,x_3\}}$,

\begin{equation}\label{E5-1-5a}
F(t,\theta) =
\begin{pmatrix}
  e^{f(t)} \tanh(t) \cos \theta \\[6pt]
  e^{f(t)} \tanh(t) \sin \theta \\[6pt]
  \dfrac{e^{f(t)}}{\cosh(t)} \\
\end{pmatrix}.
\end{equation}\bigskip

We will use the notation

\begin{equation}\label{E5-1-5b}
F(t,\theta) =
\begin{pmatrix}
e^{f(t)} \tanh(t) \, \omega_{\theta} \\[6pt]
\dfrac{e^{f(t)}}{\cosh(t)} \\
\end{pmatrix}
\end{equation} \bigskip

where $\omega_{\theta} = ( \cos \theta , \sin \theta )$ for short,
and we denote $(- \sin \theta , \cos \theta )$ by
$\dot{\omega}_{\theta}$.\bigskip

The metric induced on $M$ from the immersion $F$ is given by the
matrix

\begin{equation}\label{E5-1-6}
G_F(t,\theta) = \begin{pmatrix}
  1 + \cosh^2(t) f_t^2(t) & 0 \\[6pt]
  0 & \sinh^2(t) \\
\end{pmatrix},
\end{equation}

where $f_t$ denotes the derivative of the function $f$ with respect
to the variable $t$. \bigskip

The unit normal vector $N_F$ to the immersion is given by

\begin{equation}\label{E5-1-7}
N_F(t,\theta) = \frac{e^{f(t)}}{\big( 1 + \cosh^2(t) f_t^2(t)
\big)^{1/2}}
\begin{pmatrix}
  -(\frac{f_t(t)}{\cosh(t)} - \frac{\sinh(t)}{\cosh^2(t)})
  \, \omega_{\theta} \\[6pt]
  f_t(t) \tanh(t) + \frac{1}{\cosh^2(t)} \\
\end{pmatrix}.
\end{equation}\bigskip

The principal curvatures of the surface $M$ with respect to $N_F$
are given by

\begin{equation}\label{E5-1-8}
\left\{%
\begin{array}{lll}
k_p(t) &=& \dfrac{f_{tt}(t) \cosh(t) + 2 f_t(t) \sinh(t) + f_t^3(t)
\cosh^2(t) \sinh(t)}{\big( 1 + \cosh^2(t) f_t^2(t) \big)^{3/2}}\, ,\\
&& \\
k_n(t) &=& \dfrac{f_t(t) \cosh^2(t)}{\sinh(t)\big( 1 + \cosh^2(t)
f_t^2(t) \big)^{1/2}}\, ,\\
\end{array}%
\right.
\end{equation} \bigskip

where $k_p$ is the curvature of the generating curve in the
hyperbolic plane $\HH^2$ (see for example \cite{ST08a}). \bigskip

\begin{notes}
\bigskip \como{Notes}

See Notes [090223].

\comf{Notes} \bigskip
\end{notes}\bigskip

Taking these computations into account, the mean curvature of the
rotation surface $M$ is given by

\begin{equation}\label{E5-1-9}
H(t) \sinh(2t) = \frac{d}{dt} \; \frac{f_t(t) \sinh(t)
\cosh^2(t)}{\big( 1 + \cosh^2(t) f_t^2(t) \big)^{1/2}}.
\end{equation} \bigskip

When $H$ is assumed to be constant, Equation (\ref{E5-1-9}) provides
a first integral for the generating curves of rotation surfaces with
constant mean curvature $H$ in the plane $\HH^2_{\{u,v\}}$. These
generating curves come in a family $C_{H,a}$ and will be called
$H$-catenaries. The corresponding surfaces $\cC_{H,a}$ will be
called $H$-catenoids. They depend on a real parameter $a$. \bigskip

More precisely, we will consider three cases, depending on the value
of the mean curvature, $H=0$, $H=1$ and $0 < H < 1$. As a matter of
fact, we could consider the cases $0 \le H < 1$ and $H=1$, but the
case $H=0$ is of particular importance. \bigskip

%5-h3-2rs1.tex

We begin by general considerations. \bigskip

\subsubsection{General computations} \bigskip

\begin{figure}[h!tb]
\begin{center}
    \includegraphics[scale=0.5, angle=0]{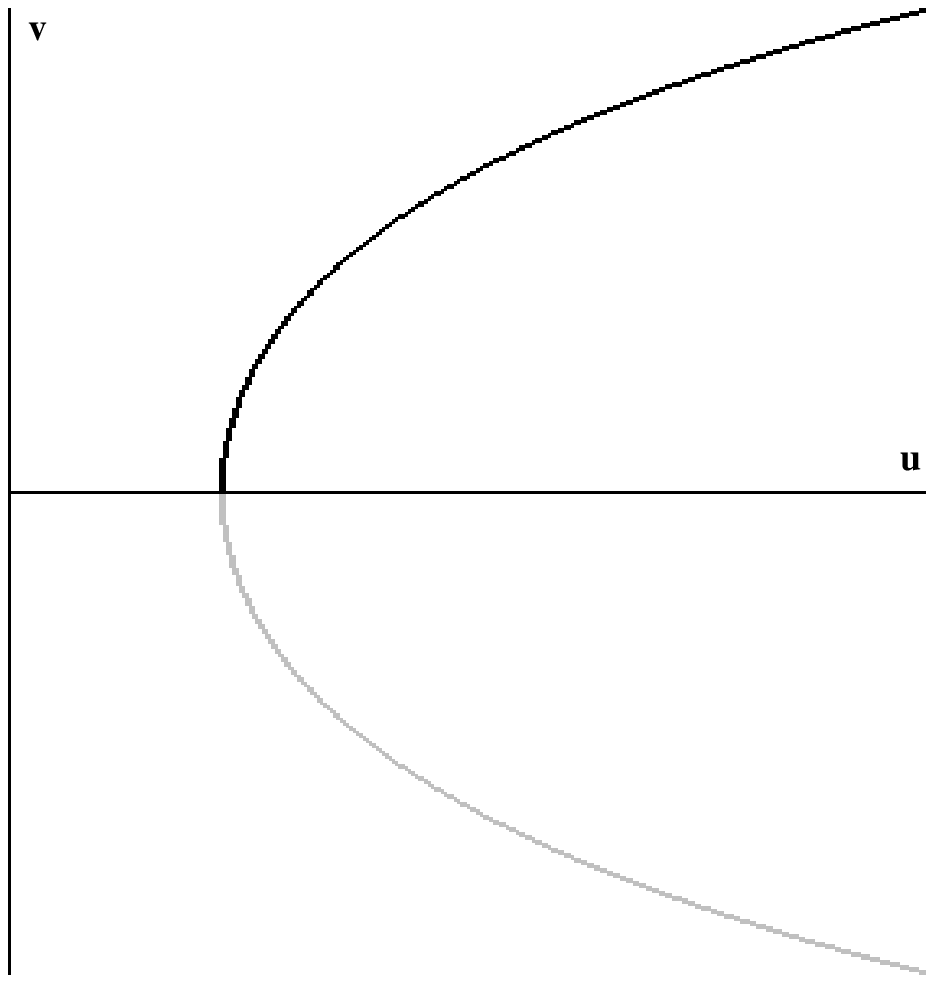}
    \caption[Curve in $\HH^2$]{Curve in $\HH^2$}
    \label{F5-3-1}
\end{center}
\end{figure}

Consider a graph $G$, $\varphi(a,t) = (t, \lambda(a,t))$ in the
plane $\HH^2_{\{u,v\}}$, see Figure \ref{F5-3-1}. Assume that the
curve extends by symmetry with respect to the $u$-axis as a smooth
curve and that the extended curve admits an arc-length
parametrization of the form

\begin{equation}\label{E5-3-1}
\Phi(a,s) = \big( y(a,s), \Lambda(a,s) \big)
\end{equation}

where $y(a,s)$ is a smooth even function of $s$ and $\Lambda(a,s)$ a
smooth odd function of $s$, such that $\Lambda(a,s) := \lambda \big(
a, y(a,s) \big)$ for $s\ge 0$.
\bigskip

\bigskip

The corresponding rotation surfaces in $\HH^3$ are given by the
parametrizations

\begin{equation}\label{E5-3-2}
\begin{array}{lll}
Y(a,s,\theta) &=&
\begin{pmatrix}
  e^{\Lambda(a,s)} \tanh(y(a,s)) \, \omega_{\theta} \\[6pt]
  \dfrac{e^{\Lambda(a,s)}}{\cosh(y(a,s))} \\
\end{pmatrix}.\\
\end{array}
\end{equation}

The parameter $s$ is the arc-length parameter along the generating
curve if and only if the following identity holds,

\begin{equation}\label{E5-3-3}
1 \equiv  y_s^2 + \cosh^2(y) \Lambda_s^2 = y_s^2 \big( 1 +
\cosh^2(y) \lambda_t^2(a,y) \big),
\end{equation}

where $y$ and $\Lambda$ stand respectively for $y(a,s)$ and
$\Lambda(a,s)$ and where subscripts indicate differentiation.

\bigskip

According to (\ref{E5-1-7}), the unit normal vectors along the
immersions are given by the formula

\begin{equation}\label{E5-3-4}
N_Y(a,s,\theta) = e^{\Lambda(a,s)}
\begin{pmatrix}
  - (\frac{\Lambda_s}{\cosh(y)} - y_s \frac{\sinh(y)}{\cosh^2(y)})
  \, \omega_{\theta} \\[6pt]
  \Lambda_s \tanh(y) + \frac{y_s}{\cosh^2(y)} \\
\end{pmatrix},
\end{equation}

where $y$ stands for $y(a,s)$. \bigskip

\begin{pb3}
\bigskip \como{PB3}

\input{pb3-h3-1}

\comf{PB3} \bigskip
\end{pb3}

Having in mind the fact that we will work with minimal or constant
mean curvature immersions, we now define Jacobi fields on the
surface $M$ in the parametrization $Y$.\bigskip

\subsubsection{Jacobi fields}\label{SS5-jf} \bigskip

Recall that the function $y(a,s)$ is assumed to be even and that the
function $\Lambda(a,s)$ is assumed to be odd.\bigskip

\noib The Killing field associated with the hyperbolic translations
along the vertical geodesic $t \mapsto (0,0,e^t)$ in
$\HH^3_{\{x_1,x_2,x_3\}}$ is just the position vector. The
\emph{vertical Jacobi field} is the function

\begin{equation}\label{E5-3-7a}
v_Y(a,s) = g_h(Y,N_Y)
\end{equation}

given by the hyperbolic scalar product of the position vector $Y$
with the unit normal vector to the immersion at $Y$. \bigskip

\begin{proper}\label{P5-3-2a}
The vertical Jacobi field $v_Y(a,s) = g_h(Y,N_Y)$ is an odd function
of $s$ given by
\begin{equation}\label{E5-3-7b}
v_Y(a,s) = \cosh(y(a,s)) y_s(a,s).
\end{equation} \bigskip
\end{proper}\bigskip

\noib The \emph{variation Jacobi field} is defined as the hyperbolic
scalar product of the variation vector-field of the family with the
unit normal vector, $e_Y(a,s) = g_h(Y_a,N_Y)$. We have

\begin{equation}\label{E5-3-9a}
Y_a(a,s,\theta) = e^{\Lambda(a,s)}
\begin{pmatrix}
  (\Lambda_a \tanh (y) + \frac{y_a}{\cosh^2(y)}) \, \omega_{\theta} \\[6pt]
  \frac{\Lambda_a}{\cosh(y)} - y_a \frac{\sinh(y)}{\cosh^2(y)} \\
\end{pmatrix}.
\end{equation}\bigskip

\begin{proper}\label{P5-3-2}
The variation Jacobi field $e_Y(a,s)$ is an even function of $s$
given by
\begin{equation}\label{E5-3-9b}
e_Y(a,s) = g_h(Y_a,N_Y) = \cosh(y(a,s)) \big( \Lambda_a y_s -
\Lambda_s y_a \big).
\end{equation}
\end{proper}\bigskip

\begin{pb3}
\bigskip \como{PB3}

\input{pb3-h3-2}

\comf{PB3} \bigskip
\end{pb3}

\begin{pb2}
\bigskip \como{PB2}

\input{pb2-h3-1}  %Jacobi fields addendum

\comf{PB2} \bigskip
\end{pb2}

%5-h3-2rs3.tex

\vspace{1cm}

We now look into the three cases, $H=0$, $H=1$ and $0 < H < 1$.
\bigskip

\subsection{Minimal catenoids in $\HH^3$} \bigskip

\subsubsection{Basic formulas}

When $H=0$, Equation (\ref{E5-1-9}) yields the solutions curves
$C_{0,a}$ (the lower index $0$ refers to the value of $H$), for
$a\ge 0$,

\begin{equation}\label{E5-1-15}
\lambda_0(a,t) = \sinh(2a) \int_a^t \dfrac{d\tau}{\cosh(\tau) \big(
\sinh^2(2\tau) - \sinh^2(2a) \big)^{1/2}}
\end{equation}

which are defined for $t \ge a$. Notice that this parametrization
only covers a half-catenary and that we work up to a $v$-translation
in $\HH^2_{\{u,v\}}$, \ie up to a hyperbolic translation with
respect to the vertical geodesic in $\HH^2_{\{x,z\}}$.
\bigskip

The arc-length parameter along the curve is given by

\begin{equation}\label{E5-1-16}
S_0(a,t) = \int_a^t \frac{\sinh(2\tau) \, d\tau}{\big(
\cosh^2(2\tau) - \cosh^2(2a) \big) ^{1/2}}
\end{equation}

or

\begin{equation}\label{E5-1-17}
\cosh (2a) \cosh \big( 2 S_0(a,t) \big) = \cosh(2t), t \ge a.
\end{equation}

\begin{prop}\label{P5-2-3}
For $s \in \R$, define the functions $y_0(a,s)$ and $\Lambda_0(a,s)$
by the formulas

\begin{equation}\label{E5-1-18}
\left\{%
\begin{array}{lll}
y_0(a,s) &=& a + \displaystyle \int_0^s \frac{\cosh(2a)
\sinh(2t)}{\big(
\cosh^2(2a) \cosh^2(2t) - 1 \big)^{1/2}} \, dt\\
\text{and}&&\\
\Lambda_0(a,s) &=& \sqrt{2} \sinh(2a) \displaystyle \int_0^s
\frac{\big( \cosh(2a) \cosh(2t) - 1 \big)^{1/2}}{\big( \cosh^2(2a)
\cosh^2(2t) - 1 \big)}
 \, dt.\\
\end{array}%
\right.
\end{equation}

\begin{enumerate}
    \item The function $y_0$ is an even function of $s$ and $\Lambda_0$ an odd
    function of $s$.
    \item For $s \ge 0$, the function $y_0(a,\cdot)$ is the inverse function
    of the function $S_0(a,\cdot)$. In particular,
    $$\cosh\big(2 y_0(a,s)\big) = \cosh(2a) \, \cosh(2s).$$
    \item For $s \ge 0$, we have $\Lambda_0(a,s) = \lambda_0(a, y_0(a,s))$.
    \item For $s \in \R$, the functions $s \mapsto \big( y_0(a,s),
    \Lambda_0(a,s)\big)$ are arc-length parametrizations of the family of
    catenaries $C_{0,a}, a>0$.
\end{enumerate}
\end{prop}\bigskip

\pf The proof is straightforward. \qed \bigskip

For later reference, we introduce the function

\begin{equation}\label{E5-10}
J_0(a,t) = \sinh(2a) (\cosh(2a) \cosh(2t) + 1)^{-1}(\cosh(2a)
\cosh(2t) - 1)^{-1/2},
\end{equation}

so that $\Lambda_0(a,s) = \sqrt{2} \int_0^s J_0(a,t) \, dt$. We
compute $\frac{\partial J_0}{\partial a}(a,t)$ and we find,

\begin{equation}\label{E5-11}
\left\{%
\begin{array}{ll}
I_0(a,t) & = \dfrac{\partial J_0}{\partial a}(a,t) =
\dfrac{n(\cosh(2a),\cosh(2t))}{d(\cosh(2a),\cosh(2t))}, \text{~where}\\
&\\
n(A,T) & = A (3 - A^2) T^2 + (A^2-1)T - 2 A,\\[6pt]
d(A,T) & = (AT + 1)^{2}(AT - 1)^{3/2}.\\
\end{array}
\right.
\end{equation}

We note that $n(A,T)$ is a polynomial of degree $2$ in $T$.

\begin{notes}
\bigskip \como{Notes}

(see Notes [090318]).

\comf{Notes} \bigskip
\end{notes}

\begin{lem}\label{L5-1}
Let $a_1 > 0$ be such that $\cosh^2(2a_1) = \frac{11+8 \sqrt{2}}{7}
\approx 3.1876$, \ie $a_1 \approx 0.5915$. For $a \ge a_1$ and for
all $t$, we have $n(cosh(2a),\cosh(2t)) \le 0$.
\end{lem}\bigskip

To the above family $C_{0,a}, a>0$ of catenaries corresponds a
family $\cC_{0,a}, a>0$ of catenoids in $\HH^3$ with the arc-length
parametrization $Y_0(a,s,\theta)$,

\begin{equation}\label{E5-12}
Y_0(a,s,\theta) =
\begin{pmatrix}
e^{\Lambda_0} \tanh(y_0)\, \omega_{\theta} \\[6pt]
e^{\Lambda_0}/\cosh(y_0) \\
\end{pmatrix}
\end{equation}

where the functions $\Lambda_0(a,s)$ and $y_0(a,s)$ are given by
Proposition \ref{P5-2-3}.
\bigskip

Catenoids in $\HH^3$ have been considered in \cite{Mo81, CD83} and
more recently in \cite{Seo09}. A new phenomenon has been pointed out
by these authors, namely that among the family $\cC_{0,a}$ of
catenoids in $\HH^3$, there are stable and index one catenoids. We
now give a precise analysis of this phenomenon and we also consider
Lindel\"{o}f's property for catenoids in $\HH^3$.
\bigskip

\subsubsection{Jacobi fields on $\cC_{0,a}$} \bigskip

According to (\ref{E5-3-4}), the unit normal $N_0(a,s,\theta)$ on
$\cC_{0,a}$ is given by

\begin{equation}\label{E5-13}
N_0(a,s,\theta) = \dfrac{e^{\Lambda_0 }}{\cosh(y_0)}\,
\begin{pmatrix}
-n_1 \, \omega_{\theta}\\[6pt]
n_2\\
\end{pmatrix}
\end{equation}

where

$$
\left\{%
\begin{array}{ll}
n_1(a,s) &= \Lambda_{0,s} - y_{0,s} \tanh(y_0), ~y_{0,s} =
\frac{\partial y_0}{\partial s},\\[8pt]
n_2(a,s) &= \Lambda_{0,s} \sinh(y_0) + y_{0,s} / \cosh(y_0).\\
\end{array}%
\right.
$$
\bigskip

Applying the formulas (\ref{E5-3-7b}) and (\ref{E5-3-9b}) of Section
\ref{SS5-jf}, we have the expressions for the vertical and variation
Jacobi fields on $\cC_{0,a}$.
\bigskip

The \emph{variation Jacobi field} $e_0(a,s)$ is given by

\begin{equation}\label{E5-14}
e_0(a,s) = - g_h(Y_{0,a}(a,s,\theta), N_0(a,s,\theta)) = -
\cosh(y_0) \big(\Lambda_{0,a}y_{0,s} - \Lambda_{0,s} y_{0,a} \big).
\end{equation}

We obtain

\begin{equation}\label{E5-16}
\begin{array}{ll}
e_0(a,s) =& \dfrac{\sinh^2(2a) \cosh(2s)}{\big( \cosh^2(2a) \cosh^2
(2s) - 1 \big)} ~~~ \cdots \\[6pt]
& \hphantom{xxxxxxx} \cdots ~~~ - \dfrac{\cosh(2a) \sinh(2s)}{\big(
\cosh(2a) \cosh (2s) - 1 \big)^{1/2}} \displaystyle \int_0^s
I_0(a,t) \, dt
\end{array}
\end{equation}

where $I_0(a,t)$ is defined by (\ref{E5-11}). \bigskip

The \emph{vertical Jacobi field} $v_0(a,s)$ is given by

\begin{equation}\label{E5-17}
v_0(a,s) = \sqrt{2} \gh(Y_0(a,s,\theta), N_0(a,s,\theta)) = \sqrt{2}
\cosh(y_0) y_{0,s}.
\end{equation}

It follows that

\begin{equation}\label{E5-19}
v_0(a,s) = \cosh(2a) \sinh(2s) \big( \cosh(2a) \cosh(2s) -
1\big)^{-1/2}.
\end{equation}

Let

\begin{equation}\label{E5-20}
f_0(a,s) = \sinh^2(2a) \cosh(2s) \big( \cosh^2(2a) \cosh^2(2s) -
1\big)^{-1}
\end{equation}

an even function of $s$ which goes to $0$ at infinity. In view of
Equations (\ref{E5-16}), (\ref{E5-19}) and (\ref{E5-20}), we have

\begin{equation}\label{E5-21}
e_0(a,s) = f_0(a,s) - v_0(a,s) \int_0^s I_0(a,t) \, dt.
\end{equation}

Observe that the integral

\begin{equation}\label{E5-21a}
    E_0(a) := \int_0^{\infty} I_0(a,t) \, dt
\end{equation}

exists for all values of $a$. \bigskip

\subsubsection{Stable domains on $\cC_{0,a}$} \bigskip

We can now investigate the stability properties of the catenoids
$\cC_{0,a}$ in $\HH^3$. \bigskip

\begin{lem}\label{L5-3-s1}
The half-catenoids
\begin{equation}\label{E5-23}
\cD_{0,a,\pm} = Y_0(a,\R_{\pm},[0,2\pi])
\end{equation}
are weakly stable. It follows from this property that a Jacobi field
$w(a,s)$ which only depends on the radial variable $s$ on
$\cC_{0,a}$ can have at most one zero on $\Rb_{+}$ and on $\Rb_{-}$.
\end{lem}\bigskip

\pf Use Property \ref{P2-4} and the fact that $v_0(a,s)$ is a Jacobi
field which only vanishes at $s=0$. \qed \bigskip

\begin{lem}\label{L5-3-s2}
The half-catenoids $Y_0(a,\R,]\varphi,\varphi +\pi[)$ are weakly
stable. Negative eigenvalues of the Jacobi operator $J_{\cC_{0,a}}$
on domains of revolution are necessarily associated with
eigenfunctions depending only on the parameter $s$. The catenoids
$\cC_{0,a}$ have at most index $1$.
\end{lem}\bigskip

\pf The fact that the index of $\cC_a$ is at most $1$ has been
proved by \cite{Seo09} using the same method as in \cite{TZ07}.
Alternatively, one could use Jacobi fields associated to geodesics
orthogonal to the axis of the catenoids. \qed \bigskip

\begin{notes}
\bigskip \como{Notes}

See Ricardo's notes \texttt{090509-rsa-lindeloef-killing.tex}.

\comf{Notes} \bigskip
\end{notes}

We can now state the main theorem of this section. Recall that the
number $E_0(a)$ is defined by  (\ref{E5-21a}) and that the Jacobi
fields $v_0(a,s)$ and $e_0(a,s)$ are given respectively by
(\ref{E5-19}) and (\ref{E5-16}), with the relation (\ref{E5-21}).
\bigskip

\begin{thm}\label{T5-1}
Let $\cC_{0,a}$ be the family of catenoids in $\HH^3$ given by
(\ref{E5-12}).
\begin{enumerate}
    \item The index of the catenoid $\cC_{0,a}$ depends on the value of the
    integral $E_0(a)$ defined by (\ref{E5-21a}). More precisely, if $E_0(a) \le  0$
    then the catenoid $\cC_{0,a}$ is stable, if $E_0(a) > 0$, then the catenoid
    $\cC_{0,a}$ has index $1$.
    \item When $\cC_{0,a}$ has index $1$, there exist $0 < z(a)$ such that
    $$\cD_{0,z(a)} = Y_0(a,]-z(a),z(a)[,[0,2\pi])$$
    is a maximal weakly stable domain.
    \item When $\cC_{0,a}$ has index $1$, there exist $0 < \ell(a) <
    z(a)$ such that
    $$\cD_{0,\ell(a)} = Y_0(a,]- \ell(a), \infty[,[0,2\pi])$$
    is a maximal weakly stable rotation invariant domain.
    \item The catenoids $\cC_{0,a}$ do not satisfy Lindel\"{o}f's property.
    \item There exist two numbers $0 < a_2 < a_1$ such that for all $a > a_1$,
    the catenoids $\cC_{0,a}$ are stable, and for all $a < a_2$, the
    catenoids $\cC_{0,a}$ have index $1$.
\end{enumerate}
\end{thm}\bigskip

\pf \bigskip

\emph{Assertion 1.} As stated in Lemma \ref{L5-3-s1}, the function
$e_0(a,s)$ can have at most one zero on $]0, \infty[$ and at most
one zero on $]- \infty,0[$. Observe that the function $e_0(a,s)$ is
even and that $e_0(a,0)=1$. To determine whether $e_0$ has a zero,
it suffices to look at its behaviour at infinity. If $E_0(a) > 0$,
the function $e_0(a,s)$ tends to $- \infty$ at infinity so that it
has exactly two symmetric zeroes in $\R$. This implies that the
index of $\cC_{0,a}$ is at least $1$. Using Lemma \ref{L5-3-s2}, we
conclude that $\cC_{0,a}$ has index $1$. If $E_0(a) < 0$, the
function $e_0(a,s)$ tends to $+ \infty$ at infinity so that it is
always positive and the catenoid $\cC_{0,a}$ is stable. Assume now
that $E_0(a)=0$. We then have the relation
$$e_0(a,s) = f_0(a,s) + v_0(a,s) \int_s^{\infty}I_0(a,t) \, dt.$$
Using Equation (\ref{E5-11}), we see that $I_0(a,t)$ is positive for
$t$ large enough provided that $\cosh^2(2a) \le 3$. In that case, it
follows that $e_0(a,s)$ is positive at infinity and hence that
$\cC_{0,a}$ is stable. If $E(a)=0$ and $\cosh^2(2a) > 3$, we need to
look at the behaviour of $e_0(a,s)$ at infinity more precisely. When
$s$ tends to $+ \infty$, we have
$$f_0(a,s) \sim 2 \tanh^2(2a) e^{-2s}, ~~v_0(a,s) \sim
\sqrt{\frac{\cosh(2a)}{2}}e^s, \text{~~~ and}$$

$$\int_s^{\infty} I_0(a,t) \, dt \sim \frac{2^{3/2}(3 - \cosh^2(2a))}
{3\cosh^{5/2}(2a)}e^{-3s}.$$

It follows that $e_0(a,s) \sim \frac{4}{3}e^{-2s}$ is positive at
infinity and hence that $\cC_{0,a}$ is stable. This proves Assertion
1. \bigskip

\emph{Assertion 2.} Saying that $\cC_{0,a}$ has index $1$ is
equivalent to saying the $E_0(a) > 0$ and hence that $e_0$ as two
symmetric zeroes. This proves Assetion 2. \bigskip

\emph{Assertion 3.} Given any $\alpha  > 0$, we introduce the Jacobi
field $e_0(a,\alpha,s)$,

\begin{equation}\label{E5-24}
e_0(a,\alpha,s) = v_0(a,\alpha) e_0(a,s) + e_0(a,\alpha) v_0(a,s).
\end{equation}

This Jacobi field vanishes at $s=-\alpha <0$ so that it cannot
vanish elsewhere in $]- \infty,0[$ and can at most vanish once in
$]0, \infty[$. Using Equations (\ref{E5-24}) and (\ref{E5-21}), we
can write

\begin{equation}\label{E5-25}
e_0(a,\alpha,s) = v_0(a,\alpha) f_0(a,s) + v_0(a,s) \big[
e_0(a,\alpha) - v_0(a,\alpha) \int_0^s I_0(a,t) \, dt \big].
\end{equation}

We have

$$e_0(a,\alpha,-\alpha) = 0 \text{~and~} e_0(a,\alpha,0) = v_0(a,\alpha) > 0$$

so that $e_0(a,\alpha, \cdot)$ vanishes in $]0, \infty[$ if and only
if $e_0(a,\alpha) - v_0(a,\alpha)E_0(a) < 0$ (recall that $E_0(a) =
\int_0^{\infty} I_0(a,t) \, dt$). \bigskip

If $\cC_{0,a}$ is stable, then clearly $e_0(a,\alpha, \cdot)$ cannot
vanish twice in $\R$. \bigskip

Assume that $\cC_{0,a}$ has index $1$ or, equivalently, that $E_0(a)
> 0$. In that case, $e_0(a, \cdot)$ has exactly one positive zero
$z(a)$. \bigskip

~~$\rhd$ For $\alpha  > z(a)$, $e_0(a,\alpha) < 0$ so that
$e_0(a,\alpha) - v_0(a,\alpha) E_0(a) < 0$ and $e_0(a,\alpha,
\cdot)$ has a positive zero $\beta$ (which must satisfy $\beta <
z(a)$).
\bigskip

~~$\rhd$ For $\alpha = z(a)$, $e_0(a,\alpha,s) =
v_0(a,\alpha)e_0(a,s)$ has two zeroes $\pm z(a)$. \bigskip

~~$\rhd$ For $0 < \alpha < z(a)$, we can argue as follows. Consider
the Jacobi field $w(a,t) = e_0(a,t) - E_0(a) v_0(a,t)$.  At $t=0$,
we have $w(a,0)=1$ and at $t=z(a)$,  we have $w(a,z(a)) < 0$ because
$e_0(a,z(a))=0$, $E_0(a) > 0$ and $v_0(a,z(a)) > 0$. It follows that
$w(a,t)$ has a unique zero in $]0,z(a)[$ and hence that there exists
a value $\ell(a) > 0$ such that

$$\cD_{0,\ell(a)} = Y_0(a, ]-\ell(a), \infty[,[0,2\pi])$$

is a maximal weakly stable rotation invariant domain. This proves
Assertion 3. \bigskip

\emph{Assertion 4.} This follows immediately from the previous
assertion. \bigskip

\emph{Assertion 5.} The first part of the Assertion follows from
Lemma \ref{L5-1} which implies that $e(a,s)$ never vanishes when $a
> a_1$. To prove the second part of Assertion 3, we can either use the fact
that $E_0(a)$ tends to $+ \infty$ when $a$ tends to zero from above
or use the criteria given in \cite{CD83} (Corollary 5.13, p. 708) or
\cite{Seo09} (Corollary 4.2), see Section \ref{SS-nc}. \qed
\bigskip

\begin{notes}
\bigskip \como{Notes}

See Ricardo's notes \texttt{catenoids-h3-v8.pdf} version 8, Theorem
1, page 9.

\comf{Notes} \bigskip
\end{notes}

We have the following geometric interpretation of Theorem \ref{T5-1}
\bigskip

\begin{figure}[htbp]
\begin{center}
\begin{minipage}[c]{7cm}
    \includegraphics[width=7cm]{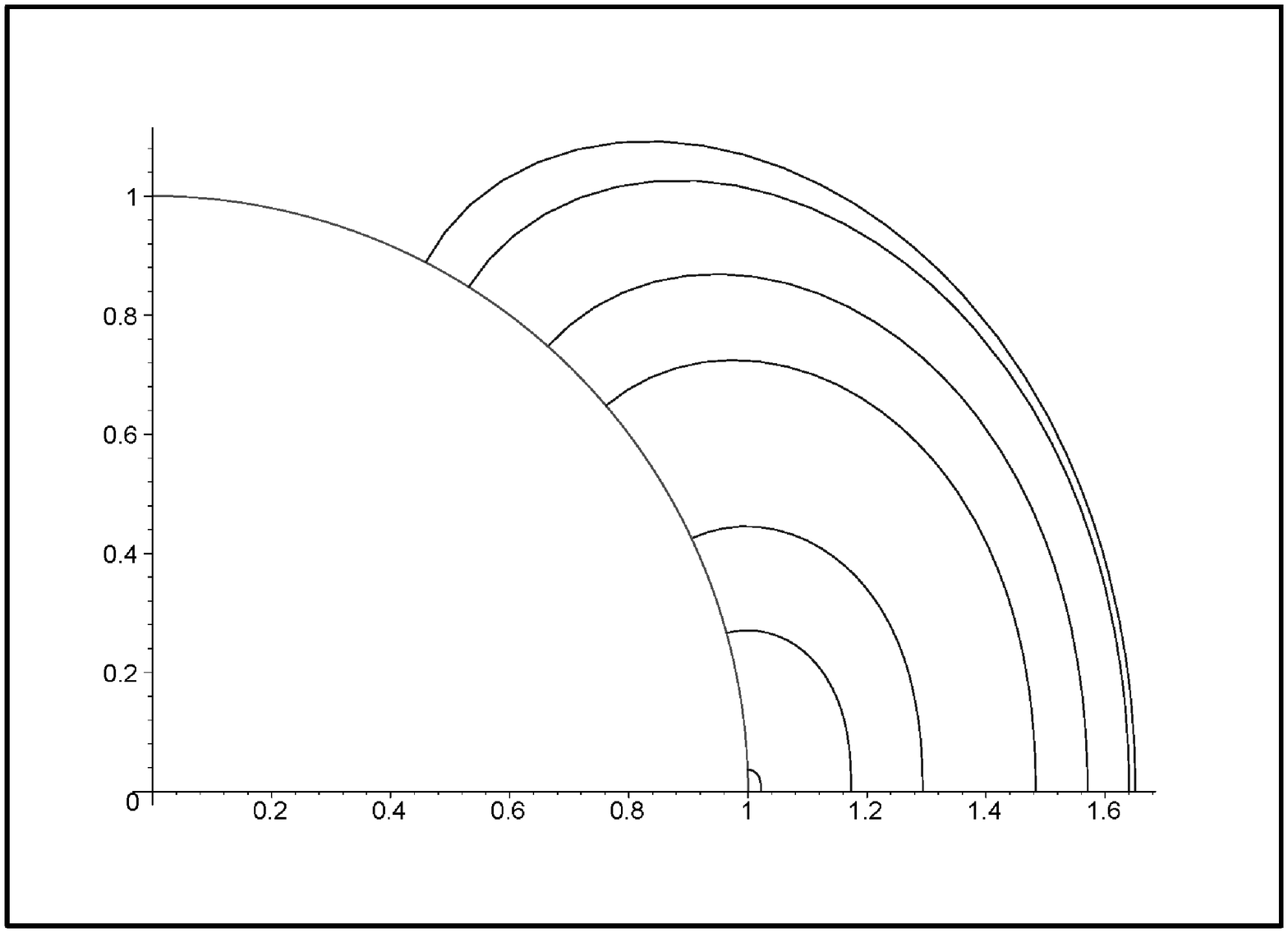}
    \caption[Foliating]{Foliating}
    \label{F5-h3-fol}
\end{minipage}\hfill
\begin{minipage}[c]{7cm}
    \includegraphics[width=7cm]{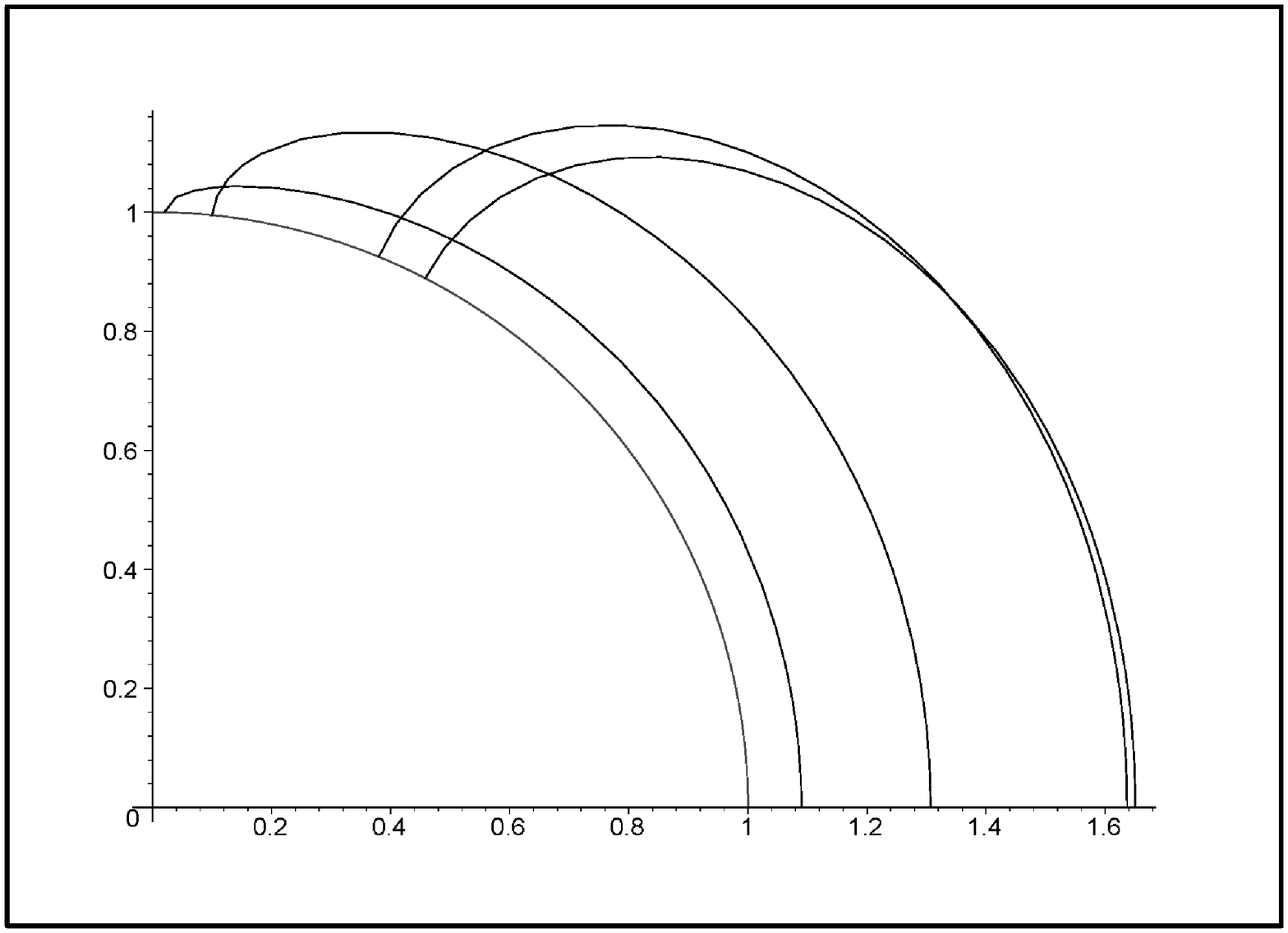}
    \caption[Intersecting]{Intersecting}
    \label{F5-h3-int}
\end{minipage}
\end{center}
\end{figure} \bigskip

\begin{prop}\label{P5-1}
We have the following geometric interpretation.
\begin{itemize}
    \item Let $\cS$ be an open interval on which $E_0 < 0$ (hence the catenoid
    $\cC_{0,a}$ is stable for all $a \in \cS$). For $a \in \cS$,
    the catenaries $C_{0,a}$ locally foliate the hyperbolic plane
    $\HH^2_{\{u,v\}}$.
    \item Let $\cU$ be an open interval on which $E_0 > 0$ (hence the catenoid
    $\cC_{0,a}$ has index $1$ for all $a \in \cU$). For $a, b \in
    \cU$, the catenaries $C_{0,a}$ and $C_{0,b}$ in $\HH^2_{\{u,v\}}$
    intersect exactly at two points. Furthermore, the family $\{C_{0,a}\}_{a
    \in \cU}$ has an envelope. Furthermore, the points at which
    $\cC_{0,a}$ touches the envelope correspond to the maximal stable
    domain $\cD_{0,z(a)}$.
\end{itemize}
\end{prop}\bigskip

\pf

Define the $v$-height function of the catenoid $\cC_{0,a}$ by

\begin{equation}\label{E-30a}
V_0(a) = \lim_{t \to \infty} \lambda_0(a,t) = \lim_{s \to \infty}
\Lambda_0(a,s).
\end{equation}

\begin{lem}\label{L5-geom}
Let $a_2 > a_1 > 0$ be two values of the parameter $a$. The
catenaries $C_{0,a_1}$ and $C_{0,a_2}$ intersect at most at two
symmetric points and they do so if and only if $V_0(a_2) >
V_0(a_1)$.
\end{lem}\bigskip

\pf To prove the Lemma, consider the difference $w(t) :=
\lambda_0(a_2,t) - \lambda_0(a_1,t)$ for $t \ge a_2 > a_1$. A
straightforward computation shows that this function increases from
the negative value $-\lambda_0(a_1,a_2)$ (achieved for $t=a_2$) to
$V_0(a_2) - V_0(a_1)$ (the limit at $t=\infty$). It follows that $w$
has at most one zero and does so if and only if $V_0(a_2) - V_0(a_1)
> 0$. \bigskip

The Proposition follows from the fact that $V_0(a) = \sqrt{2}
\int_0^{\infty}J_0(a,t) \, dt$ and that $V'_0(a) = \sqrt{2}E_0(a)$
where $E_0(a)$ is defined by (\ref{E5-21a}). \qed \bigskip

\begin{figure}[htbp]
%\begin{center}
\begin{minipage}[c]{7cm}
    \includegraphics[width=7cm]{Famindice1semEnvel}
    \caption[Intersection]{Intersecting}
    \label{F5-h3-int2}
\end{minipage}\hfill
\begin{minipage}[c]{7cm}
    \includegraphics[width=7cm]{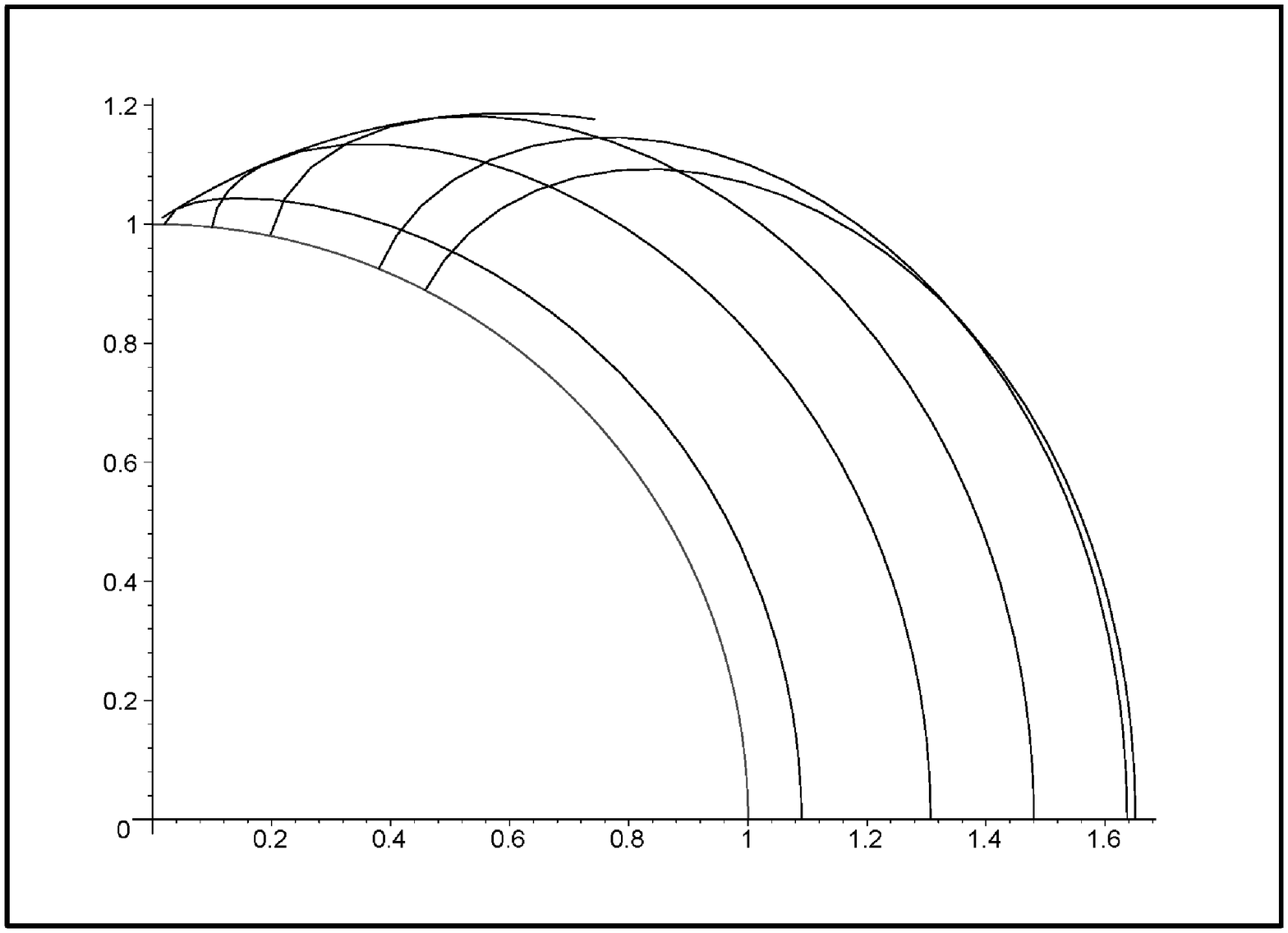}
    \caption[With envelope]{With envelope}
    \label{F5-h3-int3}
\end{minipage}\hfill
%\end{center}
\end{figure} \bigskip

\textbf{Observation.}~ One can also define the $x$-height function
of the catenoid $\cC_{0,a}$ by

\begin{equation}\label{E5-30}
X_0(a) = \lim_{s \mapsto \infty} e^{\Lambda_0(a,s)} \tanh(y_0(a,s))
= e^{\sqrt{2}\int_0^{\infty}J_0(a,t) \, dt},
\end{equation}

where $J_0(a,t)$ is defined by (\ref{E5-10}).\bigskip

The critical points of $X_0(a)$ correspond to the zeroes of the
function $E_0(a)$. \bigskip

\begin{notes}
\bigskip \como{Notes}

See [090318] for more details.

\comf{Notes} \bigskip
\end{notes}

\subsubsection{Numerical computations}\label{SS-nc} \bigskip

\textbf{Remarks.}\\

\begin{figure}[htbp]
\begin{center}
\begin{minipage}[c]{6cm}
    \includegraphics[width=4.0cm,angle=-90]{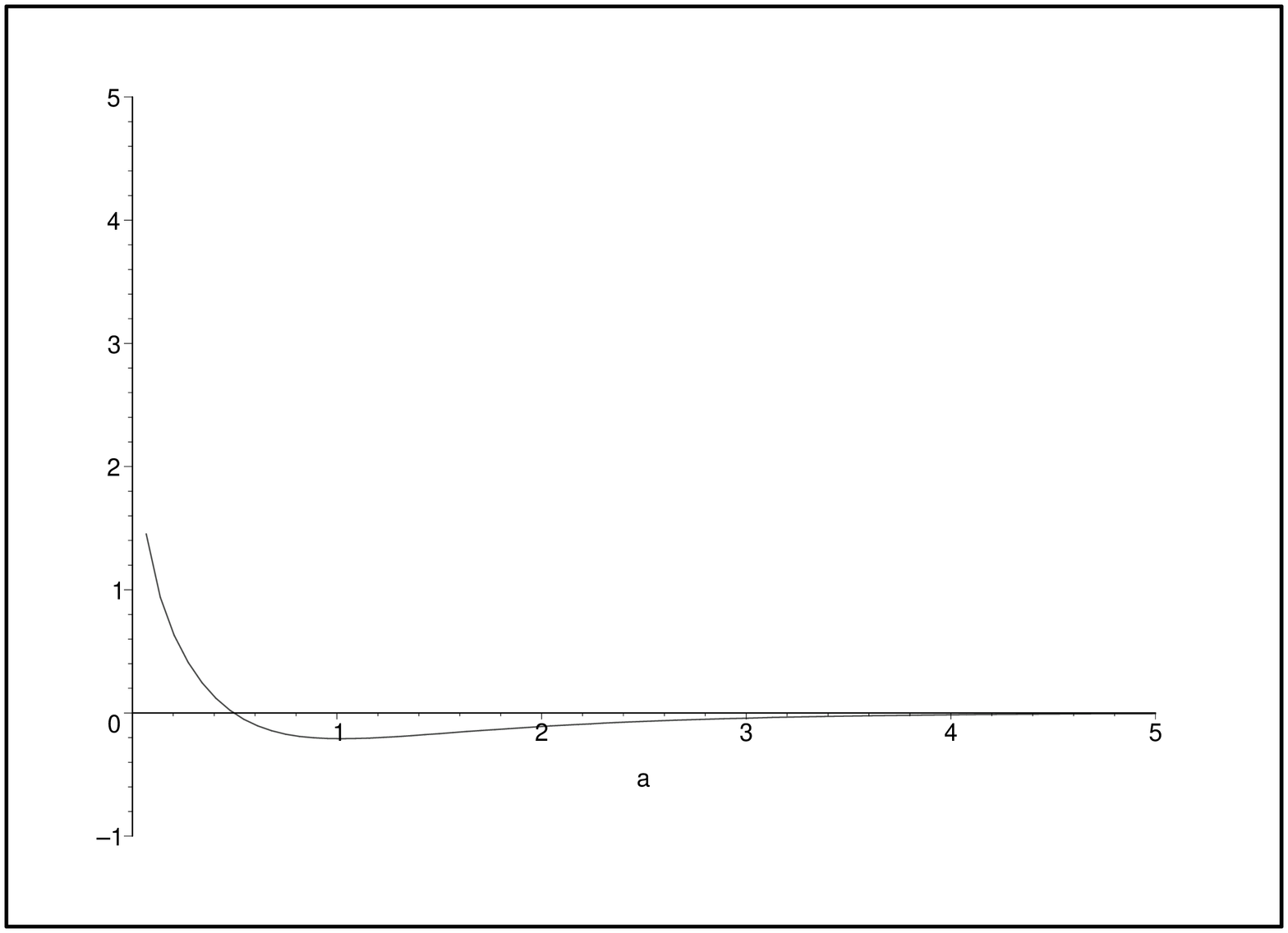}
    \caption[Graph $E_0(a)$]{Graph $E_0(a)$}
    \label{F5-h3-E}
\end{minipage}\hspace{1cm}%\hfill
\begin{minipage}[c]{6cm}
    \includegraphics[width=4.0cm,angle=-90]{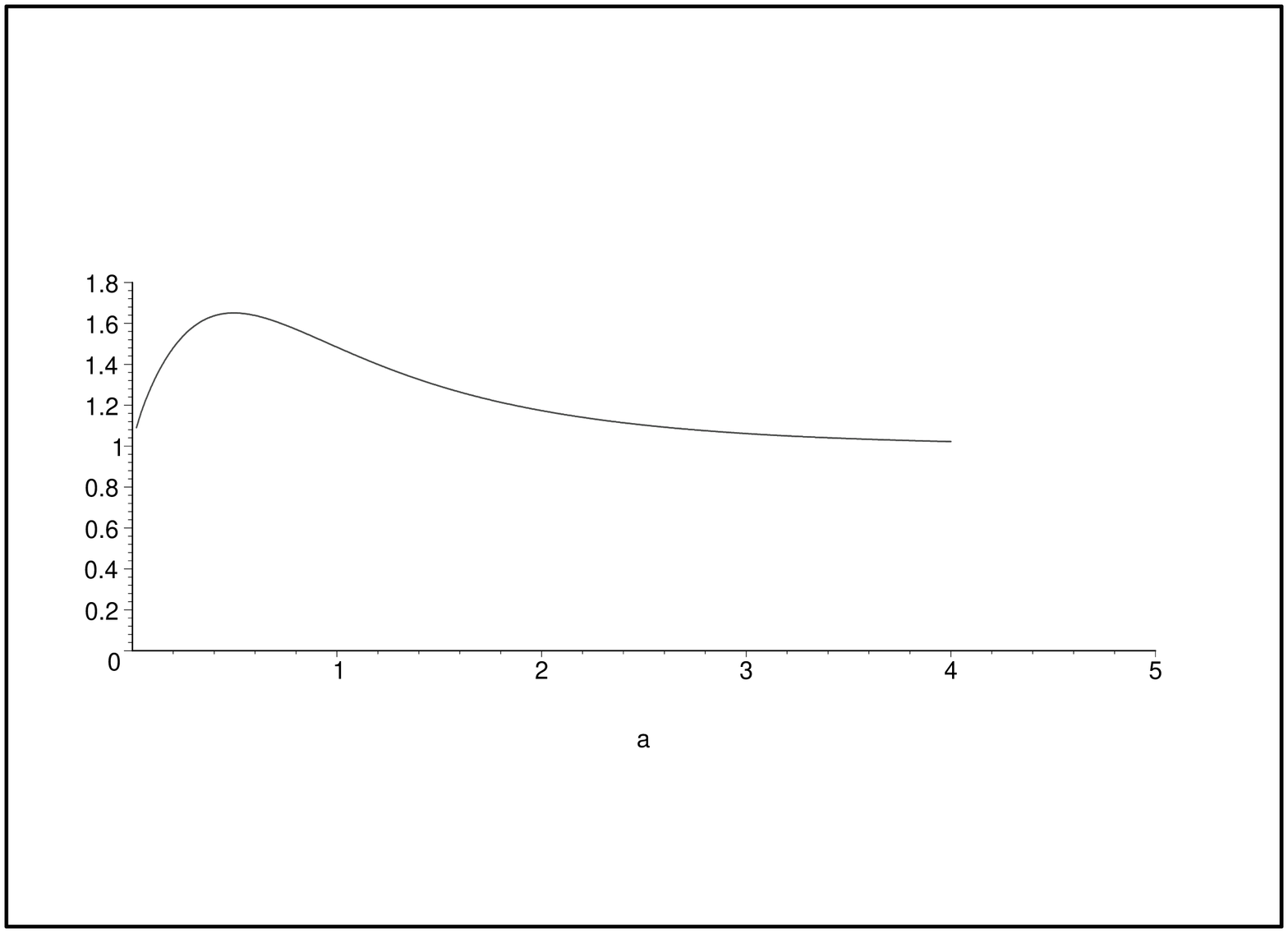}
    \caption[Graph $X_0(a)$]{Graph $X_0(a)$}
    \label{F5-h3-X}
\end{minipage}\hfill
\end{center}
\end{figure}

{1.~} The graph (\textsc{maple} plot) of the function $a \mapsto
E(a)$, for $a > 0$ (see Figure \ref{F5-h3-E}) shows that there
exists some $a_0 \approx 0,4955 \cdots $ such that

\begin{itemize}
    \item for $a \ge a_0$, the catenoids $\cC_{0,a}$ are stable and the
    corresponding catenaries locally foliate the hyperbolic plane,
    \item for $a < a_0$, the catenoids $\cC_{0,a}$ have index $1$,
    \item the function $E_0(a)$ has a unique zero $a_0$ which is the
    unique critical point of $x$-height function $X_0(a)$. The
    properties of the family of catenoids change at the
    point $a_0$, from an intersecting family to a foliating family.
\end{itemize}\bigskip

\bigskip

{2.~} The family of minimal catenoids in $\HH^3$ has been described
for the first time by H. Mori. To perform the computations, he used
the representation of $\HH^3$ as a hypersurface in the
$4$-dimensional Lorentz space. The family of catenoids is described
by a function $\phi(\alpha,s)$ (\cite{Mo81}, Theorem 1, p. 791)
which is the same as our function $\Lambda_0(a,s)$, Proposition
\ref{P5-2-3}, Equation (\ref{E5-1-18}), if we set $2\alpha =
\cosh(2a)$. According to Mori's Theorem 2 (\cite{Mo81}, p. 792), for
$\alpha \ge \frac{17}{2}$, \ie for $a \ge a_M := \mathrm{argcosh}(3)
\approx 1.7627 \cdots$, the catenoid $\cC_{0,a}$ is (globally)
stable. Mori's proof relies on the following facts.
\begin{itemize}
    \item The Jacobi operator on $\cC_{0,a}$ is given by $J = -\Delta
    + 2 - |A|^2$, where the norm of the second fundamental form
    $|A|$ can be expressed in terms of $a, s$.
    \item The Laplacian is bounded from below, $-\Delta \ge \frac{1}{4}$
    on $\cC_{0,a}$ (this follows from Cheeger's inequality).
    \item For $a \ge a_M$, we have $|A|^2 \le 2+ \frac{1}{4}$.
\end{itemize}

This method for proving stability is far from optimal. This
explains why Mori's bound $a_M \approx 1.7627 \cdots$ is worse
than our bounds $a_1 \approx 0.5915$ and $a_0 \approx 0,4955
\cdots $. \bigskip

M. do Carmo and M. Dajczer proved that some catenoids in the family
$\cC_{0,a}$ are not stable. For that purpose, they proved
(\cite{CD83}, Corollary 5.13, p. 708) that a stable complete minimal
immersion, $M^n \looparrowright \HH_{n+1}$ with finite total
curvature must satisfy the inequality

$$\int_M |A_M|^2 \big( |A_M|^2 - n(n+1) \big) \, d\mu_M \le 0.$$

Taking the explicit form of $|A|$ and $d\mu_M$ on the catenoids
$\cC_{0,a}$, the left-hand side of the above inequality give a
function of $a$. Plotting this function, one sees that the
catenoids $\cC_{0,a}$ have at least index $1$ for $a \le a_{CD}
\approx 0,4668 \cdots$. Using a different criterion, K. Seo
slightly improved this bound. The bound $a_{CD} \approx 0,4668
\cdots$ is slightly less that our bound $a_0 \approx 0,4955 \cdots
$. \bigskip

~$\maltese$~ Maple \texttt{mori-et-alii.mws} to be checked !
~$\maltese$~
\bigskip

%\subsection{Surfaces with constant mean curvature $1$} \bigskip

\subsection{Catenoid cousins in $\HH^3$} \bigskip

\subsubsection{Basic formulas}

We now consider \emph{catenoid cousins}, \ie constant mean curvature
$1$ rotation hypersurfaces in $\HH^3(-1)$. In this case, the mean
curvature equation (\ref{E5-1-9}) reads

\begin{equation}\label{E5-42}
\sinh(2t) = \frac{d}{dt} \; \frac{f_t(t) \sinh(t) \cosh^2(t)}{\big(
1 + \cosh^2(t) f_t^2(t) \big)^{1/2}}.
\end{equation}

which yields

$$\frac{f_t(t) \sinh(t) \cosh^2(t)}{\big( 1 + \cosh^2(t) f_t^2(t)
\big)^{1/2}} = \frac{1}{2} \cosh(2t) -d,$$

for some constant $d \in \R$. It follows that

$$f_t^2 \cosh^2(t) \big[ \sinh^2(2t) - (\cosh(2t) -2d)^2\big] =
\big[ \cosh(2t) -2d\big]^2$$

$$f_t^2 \cosh^2(t) \big[ 4d \cosh(2t) - 1 - 4d^2\big] =
\big[ \cosh(2t) -2d\big]^2.$$

For a solution to exist, $d$ needs to be positive so that we may
assume that $2d = e^{-2a}$ for some $a \in \R$, and we get

$$2 f_t^2 \cosh^2(t) e^{-2a} \big( \cosh(2t) - \cosh(2a) \big) =
\big( \cosh(2t) - e^{-2a} \big)^2.$$

It follows from (\ref{E5-42}) that

\begin{equation}\label{E5-44}
f_t = \dfrac{e^a \big( \cosh(2t) - e^{-2a} \big)}{\sqrt{2} \cosh(t)
\sqrt{\cosh(2t) - \cosh(2a)}}, ~~ t \ge |a|.
\end{equation} \bigskip

$\blacktriangleright$~ We now limit ourselves to the embedded case
and assume that $a>0$.
\bigskip

Equation (\ref{E5-44}) yields embedded catenary cousins $\{ C_{1,a},
a > 0\}$, given by

\begin{equation}\label{E5-46}
\lambda_1(a,t) = \int_a^t \dfrac{e^a \big( \cosh(2\tau) - e^{-2a}
\big)}{\sqrt{2} \cosh(\tau) \sqrt{\cosh(2\tau) - \cosh(2a)}} \,
d\tau, \text{ ~for~ }t \ge a,
\end{equation}

where the lower index $1$ refers to $H=1$. \bigskip

Notice that the function $\lambda_1$ describes the upper halves of
catenary-like curves and that we work up to $v$-translations in
$\HH^2_{\{u,v\}}$, \ie up to hyperbolic translations along the
vertical geodesic $\gamma$ in $\HH^2_{\{x,z\}}$. \bigskip

The arc-length function along the curve $C_{1,a}$ is given by

$$S_1(a,t) = \int_a^t \big( 1 + \cosh^2(\tau) \lambda_{1,\tau}^2(a,\tau)
\big)^{1/2} \, d\tau$$

\ie

$$S_1(a,t) = \int_a^t \frac{e^a \sinh(2\tau)}{\sqrt{2}
\sqrt{\cosh (2\tau)-\cosh(2a)}} \, d\tau.$$

Finally, we arrive at

$$S_1(a,t) = \frac{e^a}{\sqrt{2}}\sqrt{ \cosh(2t) - \cosh(2a)}$$

\ie

\begin{equation}\label{E5-48}
\cosh(2t) = 2 e^{-2a} S_1^2(a,t) + \cosh(2a), ~~ t\ge a.
\end{equation}\bigskip

For $s>0$, we can define a positive function $y_1(a,s)$ by the
relation

$$\cosh\big( 2y_1(a,s) \big) = 2 e^{-2a} s^2 + \cosh(2a)$$

and we can compute the derivative of the function $s \mapsto
\lambda_1\big( a,y_1(a,s) \big)$. We obtain the formula

$$\partial_s \lambda_1\big( a,y_1(a,s) \big) = \frac{\sqrt{2}\big( 2 e^{-2a}s^2+
\sinh(2a)\big)}{(2 e^{-2a}s^2+ \cosh(2a)+1)\sqrt{2 e^{-2a}s^2+
\cosh(2a)-1}}$$

which we can write as

$$\partial_s \lambda_1\big( a,y_1(a,s) \big) = \frac{e^a \big( 2s^2+
e^{2a} \sinh(2a)\big)}{2(s^2+ e^{2a}\cosh^2(a))\sqrt{ s^2+
e^{2a}\sinh^2(a)}}.$$

We can use this formulas to define the functions $y_1(a,s)$ and
$\Lambda_1(a,s)$ over $\R$ as follows. \bigskip

\begin{prop}\label{P5-10}
For $a>0$ and $s \in \R$, define the functions $y_1(a,s)$ and
$\Lambda_1(a,s)$ by the formulas
\begin{equation}\label{E5-50}
y_1(a,s) = a + \int_0^s \frac{2 e^{-2a}t \, dt}{\sqrt{ \big( 2
e^{-2a}t^2+\cosh(2a)\big)^2-1}},
\end{equation}
and
\begin{equation}\label{E5-52}
\Lambda_1(a,s) = \int_0^s \frac{e^a \big( 2t^2+ e^{2a}
\sinh(2a)\big) \, dt}{2\big( t^2+ e^{2a}\cosh^2(a) \big) \sqrt{ t^2+
e^{2a}\sinh^2(a)}}.
\end{equation} \bigskip
\begin{enumerate}
    \item The function $y_1$ is smooth, even, and satisfies $$\cosh(2y_1(a,s)) =
    2 e^{-2a}s^2+ \cosh(2a).$$
    \item The function $\Lambda_1$ is smooth, odd, and satisfies
    $\Lambda_1(a,s) = \lambda_1(a,y_1(a,s))$ for $s \ge 0$.
    \item For $a>0$, the maps $\R \ni s \mapsto \big( y_1(a,s),
    \Lambda_1(a,s)\big) \in \HH^2_{\{u,v\}}$ are arc-length parametrizations
    of the family of embedded catenary cousins $\{C_{1,a}\}_{a>0}$ which generate the family
    $\{\cC_{1,a}\}_{a>0}$ of embedded catenoid cousins (rotation surfaces with constant
    mean curvature $1$ in $\HH^3(-1)$).
    \item The parametrization of the family $\{\cC_{1,a}\}_{a>0}$ in
    $\HH^3_{\{x_1,x_2,x_3\}}$, is given by
    \begin{equation}\label{E5-53}
    Y_1(a,s) =
    \begin{pmatrix}
    e^{\Lambda_1(a,s)} \tanh(y_1(a,s)) \, \omega_{\theta} \\[6pt]
    e^{\Lambda_1(a,s)} \frac{1}{\cosh(y_1(a,s))}\\
    \end{pmatrix}.\bigskip
    \end{equation}
\end{enumerate}
\end{prop}\bigskip

\subsubsection{Jacobi fields on $\cC_{1,a}$} \bigskip

As in Section \ref{SS5-jf}, we define the vertical and variation
Jacobi fields on $\cC_{1,a}$. \bigskip

The \emph{vertical Jacobi field} $v_1(a,s)$ on $\cC_{1,a}$ is the
scalar product of the Killing field of hyperbolic translations along
the vertical geodesic $\gamma$ with the unit normal vector to the
surface. According to formula (\ref{E5-3-7b}), we have

\begin{lem}\label{L5-jfv1}
The vertical Jacobi field $v_1$ is a smooth odd function of $s$. It
is given by
\begin{equation}\label{E5-54}
v_1(a,s) = \cosh(y_1(a,s)) y_{1,s}(a,s) = \dfrac{e^{-a}s}{\sqrt{s^2+
e^{2a}\sinh^2(a)}}
\end{equation}
and satisfies $v_1(a,0) = 0$, $v_1(a, \infty)=e^{-a}$.
\end{lem}\bigskip

The \emph{variation Jacobi field} $e_1(a,s)$ on $\cC_{1,a}$ is the
scalar product of the variation field of the family $\cC_{1,a}$ with
the unit normal vector to the surface. According to (\ref{E5-3-9b}),
we have

\begin{equation}\label{E5-58}
e_1(a,s) = \cosh\big( y_1(a,s) \big) \big( \Lambda_{1,a}y_{1,s} -
\Lambda_{1,s} y_{1,a}\big)(a,s).
\end{equation}

which we can write as

$$v_1(a,s) \Lambda_{1,a}(a,s) - \cosh(y_1(a,s)) y_{1,a}(a,s) \Lambda_{1,s}(a,s).$$

Using Proposition \ref{P5-10} and Lemma \ref{L5-jfv1}, we find the
formula

\begin{equation}\label{E5-60}
\cosh(y_1) \, \Lambda_{1,s} \, y_{1,a} = \dfrac{\sinh^2(2a) -
4e^{-4a}s^4}{4 \big( e^{-2a}s^2 + \cosh^2(a)\big)\big( e^{-2a}s^2 +
\sinh^2(a)\big)}.
\end{equation} \bigskip

By (\ref{E5-52}), we can write $\Lambda_1(a,s)$ as $\int_0^s A(a,t)
\, dt$, where the integrand $A(a,t)$ is

\begin{equation}\label{E5-62}
\left\{%
\begin{array}{lll}
A(a,t)& = & \dfrac{2 e^a t^2 + e^{3a} \sinh(2a)}{2 \big( t^2 +
e^{2a} \cosh^2(a)\big) \big( t^2 + e^{2a} \sinh^2(a) \big)^{1/2}},\\[8pt]
& = & \dfrac{A_1(a,t)}{2 A_2(a,t) A_3^{1/2}(a,t)},\\
\end{array}%
\right.
\end{equation}

where the second equality defines the functions $A_i$. \bigskip

One can now compute the derivative of $A(a,t)$ with respect to the
variable $a$. \bigskip

$$A_a(a,t) = \dfrac{A_{1,a}(a,t)}{2 A_2(a,t) A_3^{1/2}(a,t)} -
\dfrac{A_{1}(a,t)B_2(a)}{2 A_2^2(a,t) A_3^{1/2}(a,t)} -
\dfrac{A_{1}(a,t)B_3(a)}{4 A_2(a,t) A_3^{3/2}(a,t)},$$

where

$$B_2(a) = \partial_a \big( e^{2a}\cosh^2(a)\big),
~~ B_3(a) = \partial_a \big( e^{2a}\sinh^2(a)\big).$$

It follows that

$$
\begin{array}{ll} A_a(a,t) = & \dfrac{2 e^a t^2 + e^{3a}(3 \sinh(2a) + 2
\cosh(2a))} {2 A_2(a,t) A_3^{1/2}(a,t)} - \dfrac{A_{1}(a,t)B_2(a)}{2
A_2^2(a,t) A_3^{1/2}(a,t)} ~~ \cdots\\[6pt]
  & \hphantom{xxxxxx}- \dfrac{A_{1}(a,t)B_3(a)}{4 A_2(a,t)
A_3^{3/2}(a,t)},\\
\end{array}
$$

\ie

\begin{equation}\label{E5-64b}
\left\{%
\begin{array}{l}
A_a(a,t) = B(a,t) - C(a,t), \text{ ~where~ }\\
B(a,t), C(a,t) > 0, \text{ ~for~ } a>0, t\in \R,\\
B(a,t) \sim \frac{e^a}{|t|}, \text{ ~at infinity,}\\
C(a,t) = O(\frac{1}{|t|^3}), \text{ ~at infinity.}\\
\end{array}
\right.
\end{equation} \bigskip

Finally, with the above notations, we can write the variation Jacobi
field as

$$
\begin{array}{ll}
e_1(a,s) = &\;  - \, \dfrac{e^{4a} \sinh^2(a) \cosh^2(a) - s^4}{
\big( s^2 + e^{2a} \cosh^2(a)\big)\big( s^2 + e^{2a}
\sinh^2(a)\big)}
\, - \, v_1(a,s) {\displaystyle \int_0^s C(a,t) \, dt} + ~~ \cdots\\
&\\
& \hphantom{xxxxx} \cdots ~~ + v_1(a,t) {\displaystyle  \int_0^s B(a,t) \, dt} \\
\end{array}
$$

We have proved,

\begin{lem}\label{L5-jfe1}
The variation Jacobi field $e_1$ is a smooth, even function of $s$
which can be written as
\begin{equation}\label{E5-58a}
e_1(a,s) = - f_1(a,s) + v_1(a,s) {\displaystyle  \int_0^s B(a,t) \,
dt},
\end{equation}
where the function $f_1$ is a smooth, even function of $s$, such
that $f_1(a,0)=1$ and $f_1(a, \infty)$ finite. Furthermore,

$$\lim_{s \to \infty} v_1(a,s) {\displaystyle  \int_0^s B(a,t) \,
dt} = + \infty.$$
\end{lem}\bigskip

\begin{pb2}
\bigskip \como{PB2} $\maltese$~ Computations to be checked ! ~$\maltese$
\bigskip

\input{pb2-x-details-jf}

\comf{PB2}\bigskip
\end{pb2}

\subsubsection{Stable domains on the embedded catenoid cousins} \bigskip

We can now investigate the stability properties of the embedded
catenoids cousins $\{\cC_{1,a}\}_{a>0}$ in $\HH^3(-1)$. \bigskip

\begin{lem}\label{L5-3-s12}
The upper and lower halves of the embedded catenoid cousins
\begin{equation}\label{E5-64a}
\cD_{1,a,\pm} = Y_1(a,\R_{\pm},[0,2\pi]),
\end{equation}
are weakly stable. It follows from this property that a Jacobi field
$w(a,s)$ which only depends on the radial variable $s$ on
$\cC_{1,a}$, for $a>0$, can have at most one zero on $\Rb_{+}$ and
on $\Rb_{-}$.
\end{lem}\bigskip

\pf Use Property \ref{P2-4} and the fact that $v_1(a,s)$ is a Jacobi
field which only vanishes at $s=0$. \qed \bigskip

\begin{lem}\label{L5-3-s22}
The vertical halves of the catenoid cousins
$Y_1(a,\R,]\varphi,\varphi +\pi[)$ are weakly stable. Negative
eigenvalues of the Jacobi operator $J_{\cC_{1,a}}$ on domains of
revolution are necessarily associated with eigenfunctions depending
only on the parameter $s$. The embedded catenoid cousins $\cC_{1,a}$
have at most index $1$.
\end{lem}\bigskip

\pf Consider Jacobi fields associated to geodesics orthogonal to the
axis of the catenoids. \qed \bigskip

\begin{notes}
\bigskip \como{Notes}

See Ricardo's notes \texttt{090509-rsa-lindeloef-killing.tex}.

\comf{Notes} \bigskip
\end{notes}

We can now state the main theorem of this section. Recall that the
Jacobi fields $v_1(a,s)$ and $e_1(a,s)$ are given respectively by
Lemmas \ref{L5-jfv1} and \ref{L5-jfe1}. \bigskip

\begin{thm}\label{T5-11}
Let $\{ \cC_{1,a}, a>0 \}$ be the family of embedded catenoid
cousins in $\HH^3$ given by the parametrization $Y_1$, Equation
(\ref{E5-53}).
\begin{enumerate}
    \item The Jacobi field $e_1(a,s)$ has exactly one positive zero
    $z_1(a)$ and the domains
    $$\cD_{1,a,z_1(a)}=Y_1(a,]-z_1(a),z_1(a)[,[0,2\pi])$$
    are maximal weakly stable domains.
    \item For any $\alpha > 0$, there exists a $\beta(\alpha) > 0$
    such that the domains
    $$\cD_{1,a,-\alpha, \beta(\alpha)}=Y_1(a,]-\alpha ,\beta(\alpha)[,[0,2\pi])$$
    are maximal weakly stable domains.
    \item In particular, the embedded catenoid cousins $\{\cC_{1,a}\}_{a>0}$ satisfy
    Lindel\"{o}f's property: the upper and lower halves of the embedded
    catenoid cousins $\cD_{1,a,\pm}$ are maximal rotationally symmetric domains.
    \item The index of the catenoid $\cC_{1,a}$ is equal to $1$.
\end{enumerate}
\end{thm}\bigskip

\pf \bigskip

\emph{Assertion 1.} As we have seen in Lemma \ref{L5-3-s12}, the
function $e_1(a,s)$ can have at most one zero on $]0, \infty[$ and
at most one zero on $]- \infty,0[$. By Lemma \ref{L5-jfe1}, the
function $e_1(a,s)$ is even, $e_1(a,0)=-1$ and $e_1(a,\infty)=
\infty$. It follows that $e_1(a,s)$ has exactly two symmetric zeroes
in $\R$. This proves Assertion 1.
\bigskip

\emph{Assertion 2.} Given any $\alpha  > 0$, we introduce the Jacobi
field $e_1(a,\alpha,s)$,

\begin{equation}\label{E5-24a}
e_1(a,\alpha,s) = v_1(a,\alpha) e_1(a,s) + e_1(a,\alpha) v_1(a,s).
\end{equation}

This Jacobi field vanishes at $s=-\alpha <0$ so that it cannot
vanish elsewhere in $]- \infty,0[$ and can at most vanish once in
$]0, \infty[$. Using Lemma \ref{L5-jfe1}, we can write

$$e_1(a,\alpha,s) = - v_1(a,\alpha) f_1(a,s) + v_1(a,s) \big( e_1(a,\alpha) +
v_1(a,\alpha) \int_0^s B(a,t) \, dt \big).$$

It follows that $e_1(a,\alpha, -\alpha)=0$, $e_1(a,\alpha,0) < 0$
and $\lim_{s \to \infty} e_1(a,\alpha,s) = + \infty$, and hence that
$e_1(a,\alpha, \cdot)$ must vanish at least once. This proves
Assertion 2.
\bigskip

\emph{Assertion 3.} This is a consequence of Assertion 2. \bigskip

\emph{Assertion 4.} This assertion follows from Assertion 1 and from
Lemma \ref{L5-3-s22}. This has also been proved, using different
methods, by Lima and Rossman \cite{LR98}. \qed \bigskip

\subsection{Surfaces with constant mean curvature $0 \le H <
1$}\bigskip

One can also study rotation surfaces with constant mean curvature
$H$, $0 \le H < 1$ in $\HH^3(-1)$. This is similar to the case of
minimal surfaces. \bigskip

More precisely, $H$-rotation surfaces in $\HH^3(-1)$, with $0 \le
H < 1$, come in a one-parameter family $\cC_{H,a}$. For some
values of $a$ the surfaces are stable, for other values of $a$
they have index $1$. Furthermore, they do not satisfy
Lindel\"{o}f's property.  \bigskip

The computations are much more complicated but similar to the
minimal case. The functions involved depend continuously on the
parameter $H$, for $0 \le H < 1$.   \bigskip

\subsection{Higher dimensional catenoids}\bigskip

The method described in the previous sections could be applied to
study the stable domains on higher dimensional catenoids (minimal
rotation hypersurfaces or constant mean curvature $1$ rotation
hypersurfaces) in $\HH^{n+1}$.  \bigskip

\label{9-flux}
%%%\section{Appendix A}\label{A-a}
%9-flux.tex

\bigskip
\section{Appendix A}\label{A-a}
\bigskip

In this Appendix, we give a flux formula which is valid in a quite
general framework. \bigskip

Let $M^n \looparrowright (\Mh^{n+1}, \gh)$ be an isometric embedding
with mean curvature vector $\overrightarrow{H}$.  \bigskip

Given a relatively compact domain $D \subset M$, let $\nu$ denote
the unit normal to $\partial D$, pointing inwards. We denote by
$d\mu_M$ the Riemannian measure for the induced metric on $M$ and
by $d\mu_{\partial D}$ the Riemannian measure for the induced
metric on $\partial D$. \bigskip

\begin{prop}\label{P-A-1}
Given any Killing vector-field $\cK$ on $\Mh$, we have the
\emph{flux formula},
\begin{equation}\label{E-A-2}
n \int_{D} \gh(\cK, \overrightarrow{H}) \, d\mu_{M} = - \;
\int_{\partial D} \gh(\cK, \nu) \, d\mu_{\partial D}.
\end{equation}
\end{prop}\bigskip

\pf Recall that according to \cite{KN63} (p. 237 ff), the
vector-field $\cK$ is a Killing field on $\Mh$ if and only if
$\gh(\Dh_X\cK,X)=0$ for all vector-field $X$ on $\Mh$ (here $\Dh$ is
the covariant derivative associated with the Riemannian metric
$\gh$). \bigskip

Given the Killing field $\cK$, let $\omega$ be the dual $1$-form,
$\omega(\cdot) = \gh(\cK,\cdot)$ and let $\omega_M = \omega|_M$ be
the restriction of $\omega$ to $M$. \bigskip

The following formula holds,

\begin{equation}\label{E-A-4}
\delta_M \omega_M = - n \; g_M(\cK,\overrightarrow{H})
\end{equation}

where $\delta_M$ is the divergence in the induced metric on $M$.
\bigskip

Indeed, let $\{E_i\}_{1\le i\le n}$ be a local \textsc{onf} on $M$,
then

$$
\begin{array}{ll}
\delta_M \omega_M & = - \; \sum_i (D_{E_i}\omega_M)(E_i) \\
& = - \; \sum_i [E_i \cdot (\omega_M(E_i)) - \omega_M(D_{E_i}E_i)] \\
& = - \; \sum_i  [E_i \cdot \gh(\cK,E_i) + \gh(\cK,D_{E_i}E_i)] \\
& = - n \gh(\cK,\overrightarrow{H}) - \; \sum_i  \gh(\Dh_{E_i}\cK,E_i)\\
\end{array}
$$

which proves Formula (\ref{E-A-4}). We can now apply the divergence
theorem in $D$,

$$
\begin{array}{ll}
n \int_{D} \gh(\cK, \overrightarrow{H}) \, d\mu_{M} & = - \;
\int_{D} \delta_M(\omega_M) \, d\mu_M \\
& = - \int_{\partial D} \omega_M(\nu) \, d\mu_{\partial D}.
\end{array}
$$

The Proposition follows. \qed \bigskip

\begin{notes}
\como{Notes}

[090712-divergence-killing-fields-flux.pdf]

\comf{Notes}
\end{notes}\bigskip
   %Flux formula

%\newpage
\vspace{2cm}

%\def\refname{References}
%\bibliographystyle{plain}
%\bibliography{pbrsa-lind}%pbrsa-lind.bib

%\newpage
\vspace{2cm}

\begin{flushleft}
Pierre B\'{e}rard\\
Universit\'{e} Joseph Fourier\\
Institut Fourier - Math\'{e}matiques (UJF-CNRS)\\
B.P. 74\\
38402 Saint Martin d'H\`{e}res Cedex - France\\
\verb+Pierre.Berard@ujf-grenoble.fr+
\end{flushleft}\bigskip

\begin{flushleft}
Ricardo Sa Earp\\
Departamento de Matem\'{a}tica\\
Pontif{\'{i}}cia Universidade Cat\'{o}lica do Rio de Janeiro\\
22453-900 Rio de Janeiro - RJ - Brazil\\
\verb+earp@mat.puc-rio.br+
\end{flushleft}\bigskip

\vfill

\begin{flushright}
%\tiny{\texttt{\finkfile ~(\today )}}
\end{flushright}\bigskip

\end{document}